\documentclass{article}
\usepackage{arxiv}

\usepackage{hyperref}
\usepackage{bm}% allocate bold fonts
\usepackage{cite}
\usepackage{graphics,fancyhdr,graphicx,subfigure}
\usepackage{subfigure}
\usepackage{amsthm,amsmath,amsfonts,latexsym,amssymb}
\usepackage{lineno}
\usepackage{diagbox}
\usepackage{algorithm}
\PassOptionsToPackage{noend}{algpseudocode}
\usepackage{algpseudocode}
\usepackage{multirow,booktabs}
\usepackage{listings}
\usepackage[table]{xcolor}
\usepackage{lscape}
\usepackage{comment}
\usepackage{tabularx}
\usepackage{graphicx}%to add images
\usepackage{tabularx}
\usepackage{arydshln}
\usepackage{nomencl}%for adding symbols page
\usepackage{array}%to add table
\usepackage{caption}%add caption for the table and list it in list of tables
\usepackage{breqn}
\usepackage{lineno}
\usepackage{mathtools}
\usepackage{subfigure}
\usepackage{caption}
\usepackage{lmodern}
\usepackage{calligra}
\usepackage{xspace}
\usepackage[utf8]{inputenc}%skiplength
\usepackage{enumerate}%list
\usepackage[english]{babel}
\def\equationautorefname~#1\null{%
  Eq.~(#1)\null
  }
\def\subfigureautorefname~#1\null{%
  Fig.~#1\null
}

\definecolor{listinggray}{gray}{0.9}
\definecolor{lbcolor}{rgb}{0.9,0.9,0.9}
\definecolor{Darkgreen}{RGB}{0,100,0}

\title{Harnessing physics-informed operators for high-dimensional reliability analysis problems}

\author{ \hspace{1mm}Navaneeth~N.\\
	Department of Applied Mechanics\\
	Indian Institute of Technology (IIT) Delhi\\
	Hauz Khas - 110 016, New Delhi, India \\
	\texttt{navaneeth.n@am.iitd.ac.in} \\
	\And
    \hspace{1mm}Tushar\\
	Department of Applied Mechanics\\
	Indian Institute of Technology (IIT) Delhi\\
	Hauz Khas - 110 016, New Delhi, India \\
	\texttt{amz218314.iitd@gmail.com} \\
	\And
	\hspace{1mm}Souvik~Chakraborty \\
	Department of Applied Mechanics\\
        Yardi School of Artificial Intelligence \\
	Indian Institute of Technology (IIT) Delhi\\
	Hauz Khas - 110 016, New Delhi, India \\
	\texttt{souvik@am.iitd.ac.in} \\
}

\renewcommand{\shorttitle}{Harnessing physics-informed operators for high-dimensional reliability analysis problems}

\hypersetup{
pdftitle={Harnessing physics-informed operators for high-dimensional reliability analysis problems},
pdfsubject={q-bio.NC, q-bio.QM},
pdfauthor={Navaneeth N., Tushar, Souvik Chakraborty},
pdfkeywords={Operator learning, Wavelet transformations, Physics informed learning, Reliability},
}
\begin{document}
\maketitle

\begin{abstract}

Reliability analysis is a formidable task, particularly in systems with a large number of stochastic parameters. Conventional methods for quantifying reliability often rely on extensive simulations or experimental data, which can be costly and time-consuming, especially when dealing with systems governed by complex physical laws which necessitates computationally intensive numerical methods such as finite element or finite volume techniques. On the other hand, surrogate-based methods offer an efficient alternative for computing reliability by approximating the underlying model from limited data. Neural operators have recently emerged as effective surrogates for modelling physical systems governed by partial differential equations. These operators can learn solutions to PDEs for varying inputs and parameters. Here, we investigate the efficacy of the recently developed physics-informed wavelet neural operator in solving reliability analysis problems. In particular, we investigate the possibility of using physics-informed operator for solving high-dimensional reliability analysis problems, while bypassing the need for any simulation. Through four numerical examples, we illustrate that physics-informed operator can seamlessly solve high-dimensional reliability analysis problems with reasonable accuracy, while eliminating the need for running expensive simulations.

\end{abstract}


\section{Introduction}
\label{sec:intro}
Uncertainties are inherent in virtually all practical systems, affecting their performance and reliability. These uncertainties typically originate from three primary sources: model formulations, model parameters, and the inputs to the models. For example, model formulations often involve assumptions and simplifications that do not fully capture the complexity of real-world systems, introducing significant uncertainty. Similarly, model parameters, often estimated from limited data or empirical observations, can vary widely, adding to the overall uncertainty. Inputs, such as environmental conditions, are also highly variable and challenging to predict accurately. To ensure the safe operation of a system, it is imperative to quantify these uncertainties and evaluate their influence on the system’s performance, often represented as either the probability of failure or reliability. Unfortunately, this is a nontrivial undertaking and is computationally expensive due to the necessity of running repeated simulations. The challenge becomes even more pronounced when dealing with systems that have a large number of random variables. Therefore, there is a critical need to develop methods and algorithms that can accurately compute the reliability of systems with a large number of random variables using a minimal number of simulations. Such approaches would significantly enhance our ability to predict system performance, mitigate risks, and ensure the robustness and dependability of practical systems.

Among the existing methods for reliability analysis \cite{ni2020reliability,matteo2021time, de2005naroas, tao2024reliability}, the Monte Carlo simulation (MCS) \cite{thakur1978monte,rubinstein2016simulation,xiong2021fast,chang2022mc} has emerged as the foremost and most direct approach to reliability analysis. It quantifies reliability through the execution of a large number of direct simulations, drawing samples independently from the probability distribution of input variables. While the method is straightforward, it becomes computationally expensive due to the necessity of a large number of actual simulations to ensure convergence of the result. 
In response to this challenge, research has led to the development of alternative variants of the crude MCS aimed at improving computational efficiency. Notable among these are 
importance sampling \cite{au1999new,li2005curse,engelund1993benchmark}, subset simulations \cite{au2001estimation,au2014engineering,zuev2015subset}, and directional simulations \cite{ditlevsen1990general}. Importance sampling focuses on drawing samples from a distribution that emphasizes the critical regions of the input space, thus improving convergence rates. Subset simulations break down the problem into a series of conditional probabilities, also enhancing convergence. Directional simulations, on the other hand, transform the problem into a lower-dimensional space to achieve better efficiency. While these methods offer superior convergence rates compared to crude MCS, they still demand a considerable number of simulations to achieve accurate estimates of the probability of failure. In addition to sampling-based methods, analytical approximation methods such as the First-Order Reliability Method (FORM) and the Second-Order Reliability Method (SORM) are often employed. These approaches approximate the multivariate integral of the limit state function over the failure domain using Taylor series expansion and asymptotic methods. FORM and SORM offer computational efficiency but do not always guarantee convergence, especially for highly nonlinear and high-dimensional systems.

A third class of methods for reliability analysis is rooted in training efficient emulators as a surrogate to the computationally expensive model \cite{xu2021machine}. Some of the popular methods belonging to this class include Polynomial chaos expansion \cite{blatman2011adaptive,sudret2008global}, radial basis function (RBF) \cite{de2013new,li2018sequential}, ANOVA-HDMR \cite{ziehn2009gui}, Multivariate Adaptive Regression Splines(MARS) \cite{metya2017system}, Gaussian process \cite{bilionis2012multi,bilionis2013multi,tripathy2016gaussian,atkinson2018structured,atkinson2019structured, bansal2022physics,qian2021time}, Active learning Kriging Monte Carlo Simulation (AK-MCS) \cite{ling2020efficient,xiong2021fast}, support vector machines (SVM) \cite{roy2020support,ghosh2018support,cheng2021adaptive,zhou2022towards,pepper2022adaptive,chen2022support}, Artificial Neural Networks (ANN) \cite{elhewy2006reliability,navaneeth2022koopman} and hybrid surrogate models \cite{navaneeth2022surrogate,chakraborty2017efficient,chakraborty2017hybrid,wu2023adaptive} to name a few. Although these surrogate-based approaches are cost-efficient in comparison with simulation-based approaches and yield generally accurate predictions, they are purely data-driven in nature and, hence, do not comply with the underlying physics of a problem. As a consequence, these models often fail to generalize beyond the training data regime. One approach to address this challenge is using an equation discovery algorithm to determine the governing physics \cite{mathpati2023mantra}, or utilizing partially known physics in conjunction with ANNs \cite{chakraborty2023deep, chakraborty2023dpa}. However, the computational cost significantly increases in the second stage when solving the identified governing equation. Another recent attempt towards alleviating this challenge includes extending physics-informed neural networks (PINN) \cite{ramabathiran2021spinn, sirignano2018dgm, samaniego2020energy, nguyen2021parametric, chiu2022can, kharazmi2021hp, gao2021phygeonet, sharma2022accelerated} to solve reliability analysis problems \cite{chakraborty2020simulation,zhang2022simulation}. The basic idea here is to train a neural network by minimizing the residual loss computed directly from the governing equations described in the form of differential equations. However, solving problems with large number of random variables by exploiting physics-informed learning remains an open problem.

One concurrent development, in parallel to the PINN, has been the operator learning framework. Operator learning \cite{jiang2024fourier} paradigms effectively learn solutions for families of parametric PDEs by learning the mapping between input and output functions. For example, DeepoNet \cite{lu2019deeponet,lu2021learning}, considered to be the first of this kind, utilizes two networks to learn the solution operator. Graph Neural Operators (GNO) \cite{li2020multipole, xia2023maintenance}, an alternative approach, proposes learning the operator through the integral transform and kernels, where the kernels are obtained as message-passing interfaces within graph networks. Other operator learning algorithms that exploit the kernel-based architecture include Fourier Neural Operator (FNO) \cite{kovachki2021universal}, Wavelet Neural Operator (WNO) \cite{tripura2023wavelet1,tripura2023elastography, rani2024generative}, and Laplace Neural Operator (LNO) \cite{cao2023lno}. Of late, this idea has also been extended to develop the first of its kind foundation model in computational mechanics \cite{tripura2023foundational}. However, all the operator learning models discussed above are data-driven in nature, and hence, the bottleneck associated with data-driven surrogate models holds true. Therefore, the direct application of these operator learning algorithms for reliability analysis remains challenging.

In response to the above challenges, we investigate the possibility of using the recently developed physics-informed operator learning algorithm \cite{navaneeth2024physics} for solving reliability analysis problems. We hypothesize that the physics-informed operator can potentially solve high-dimensional reliability analysis problems while completely eliminating the requirement of generating data by running computationally expensive simulations. In this study, we particularly investigate the potential of using the physics-informed wavelet neural operator (PI-WNO) \cite{navaneeth2024physics} for solving high-dimensional reliability analysis problems. We investigate the applicability of PI-WNO for solving both time-independent and time-dependent reliability analysis problems. 

The remainder of the paper is organized as follows. The general problem statement is presented in Section \ref{sec:Problem statement}. Section \ref{sec: Methedology} describes the details of the methodology. Numerical examples are exemplified in the section\ref{sec: Numerical example}. Finally, the summary and concluding notes are provided in Section \ref{sec:Conclusions}.

\section{Problem statement}\label{sec:Problem statement}
In the context of reliability analysis, systems can be broadly classified into time-dependent and time-independent systems. We first consider a time-independent system characterised by an N-dimensional vector of random variables $\bm{\Lambda}$, $\bm{\Lambda}=\left(\lambda_{1}(\bm X), \lambda_{2}(\bm X),\cdots, \lambda_{N}(\bm X),\right)$: $\Omega \to \mathbb{R}^{N}$, where $\bm X$ represents the space over which the inputs and parameters are defined. We then define a mapping $\hat{J}$ that maps that N-dimensional input space to the corresponding output response. For a time-independent system having response $\hat{J}(\bm{\Lambda},\bm{X})$ and threshold $e_h$, the failure probability is quantified based on the limit state function such that $e_h - \hat{J}(\bm{\Lambda},\bm{X}) = 0 $, denotes the limiting condition. The failure domain ($\Omega^{F}$) satisfies condition: $\hat{J}(\bm{\Lambda},\bm{X})>e_h$, and is expressed as: 
\begin{equation}
\Omega^{F} \overset{\Delta}{=}\{\bm \Lambda : e_h -\hat{J}(\bm \Lambda)<0\},
\end{equation}
where the  safe region satisfies the condition: $\hat{J}(\bm \Lambda,\bm X) < e_h $. The failure probability is defined as:
\begin{equation}
P_{f}=\mathbb P(\bm{X}\in \Omega^{F})=\int_{\Omega_{\bm{\Lambda}}^{F}}dF_{\bm{\bm \Lambda}}(\bm \lambda, \bm x).
\end{equation} 

Here, $\mathbb P$ is the probability, $F_{\bm{\Lambda}}\left(\bm{\lambda}, \bm x\right)$ represents a cumulative distribution function for the probability density function $P_{\bm \Lambda}(\bm \lambda, \bm x)$ such that $F_{\bm{\Lambda}} \left(\bm{\lambda}\right)= \mathbb P \left(\bm{\Lambda} \leq \bm{\lambda}, \bm{X} \leq \bm{x} \right)$ on a given probability space, $\Omega_{\bm X}$.

In contrast to the time-independent system, for time-independent systems, output response, obtained by the mapping $\hat{J}$, depends not only on the system input and parameters ($\bm {\Lambda}$) but also varies with the time. Moreover, the threshold of the system ($e_h$) can also be a function of time. Thus, for a time-dependent system, failure probability is obtained as:
\begin{equation}
    P_{f}(t)=\mathbb P(\hat{J}(\bm \Lambda,\bm {X},t)>e_h(t)).
\end{equation}
While there are several approaches to obtaining the quantitative assessment of time-dependent failure and time-dependent reliability, in the present work, we employ the first passage failure method to estimate the time-dependent reliability. The first passage failure time, also known as First Time To Failure (FTTF), is defined as the time ($\tau$) at which the system response, $\hat{J}({\cdot})$ crosses the threshold for the first time. When the system is considered in between  time intervals with the $t_0$ and $t_s$ being the initial and final points, the first passage failure probability can be defined as follows:
\begin{equation}\label{eq:time_threshold}
  P_{f}=\mathbb P(\hat{J}(\bm \Lambda,\bm{X},\tau)<e_h(\tau), \tau\in [t_o,t_s]),
\end{equation}
The overarching objective here is to investigate the possible application of physics-informed operators for solving time-independent and time-dependent reliability analysis problems in the high dimensional system.
\section{Physics informed operator}\label{sec: Methedology}
In this section, we provide a detailed description of the physics-informed WNO and the implementation of physics-informed operator for reliability analysis. 

\subsection{Operator learning}
Operator learning aims to learn a function to function mapping from data. Consider a system governed by a parametric PDE of the following form: 
\begin{equation}\label{operator1}
    \mathcal{N}(\bm{a},\bm{u},\bm x,t) = \bm s(\bm x,t) , \,\, \text { in } D \subset \mathbb{R}^d.
\end{equation}
The PDE is defined on an $d$-dimensional domain, $D \in \mathbb{R}^{d}$ bounded by $\partial D$, with the initial conditions and boundary conditions, respectively, being of the forms:
\begin{equation}\label{genBC}
\begin{aligned}
\bm{u}\left(\bm{x}, 0\right) &=\bm{u}_0(\bm{x}) \quad \bm{x} \in D, \\
\bm{u}(\bm{x}, t) &= \bm{g}(\bm x,t), \quad \bm{x} \in \partial D,\, t \in[0, T],
\end{aligned}
\end{equation}
where $\bm x$ and $t$ represent the space and the time coordinates, respectively. Here, $a$, which belongs to the function space $\mathcal{A}$ , i.e., $\bm{a} \in \mathcal{A} $, denotes the set of input parameters, with $\bm a:D \mapsto \mathbb{R}^{da}$. $\bm u$ is the solution of the nonlinear differential operator $\mathcal{N}$  corresponding to the set of inputs $\bm a$, and source function $\bm{s}$ subjected to boundary conditions $\bm{g}$ and initial conditions $\bm{u}_0$ such that $\bm{u} \in \mathcal{U} $. The output function space comprises the solution of the parametric PDE, $u(x, t): D \mapsto \mathbb{R}$. 
For the differential operator $\mathcal{N}$ and the aforementioned inputs, and given input conditions, there exists a solution operator $\mathcal{M}: \mathcal{A} \mapsto \mathcal{U}$, which maps the input functions to the solution. The nonlinear solution operator, $\mathcal{M}$, approximated by the neural network can be expressed as: 
\begin{equation}\label{neural operator}
    \bm{u} \approx \hat{\bm{u}} = \hat{\mathcal{M}} (\bm{a},\bm{x};\bm{\theta}_{NN})
\end{equation}
where $\bm{\theta}_{NN}$  represents the trainable neural network parameters of operator
$\hat{\mathcal{M}}$. 
Here, we note the underlying assumption for the operator learning that for any $\bm{a} \in \mathcal {A}$, there exists a unique solution $\hat {\bm{u}} = \hat{\mathcal{M}}(\bm{x}, \bm{a};\bm{\theta}_{NN}) \in \mathcal{U}$. The different operator learning algorithms existing in the literature differ in how the function $\hat {\mathcal M}$ is parameterized. For example, in DeepONet, the function $\hat {\mathcal M}$ is parameterized using two networks, the trunk net and the branch net. Similarly, in kernel-based operator, the function $\hat {\mathcal M}$ is parameterized using kernel integration. In this work, we use a wavelet neural operator to represent $\hat {\mathcal M}$. Accordingly, details on the wavelet neural operator are presented next.

\subsection{Wavelet neural operator}
The Wavelet Neural Operator (WNO) is a kernel-based neural operator that aims to learn the integral operator for a family of parametric partial differential equations (PDEs). Consider a parametric PDE with the input function space denoted as \(\mathcal{A}\) and the output space denoted as \(\mathcal{U}\). The neural operator is responsible for learning the integral operator $\hat{\mathcal{M}}$ such that $u(x) \approx \hat{\mathcal{M}}(a(x),x)$, with the input-output pair $\{a(x) \in \mathcal{A}, u(x) \in \mathcal{U}\}$ defined on a smooth $d$-dimensional domain $D$. Kernel-based operator learning methods are inspired by the  Hammerstein integral equation,
\begin{equation}\label{eq:hammerstein}
    u(x) = \int_{D} k(x, \xi)f\left(\xi, u(\xi)\right) d\xi + m(x); \quad x \in D,
\end{equation}
where $m$ represents a nonlinear transformation and $k(\cdot)$ denotes the kernel of the nonlinear integral equation. Here, $f(x,\xi)$ is a given input function that is continuous and satisfies the conditions (Lipschitz criteria):
\begin{equation}
\begin{aligned}
    |f(x,\xi)| & \leq C_1|\xi|+C_2\\
    \left|f\left(x, \xi_1\right)-f\left(x, \xi_2\right)\right| &\leq C\left|\xi_1-\xi_2\right|
\end{aligned}
\end{equation}
where the constants, $C$, $C_1$ and $C_2$ are always positive with $C$, $C_1$ being smaller than the first eigen value of the kernel $K(x,\xi)$. \autoref{eq:hammerstein} is theoretically defined on an infinite-dimensional space, whereas in practice, it is defined on a finite-dimensional parameterized space obtained by discretizing the solution domain $D \in \mathbb{R}^n$. To facilitate learning complex features through subsequent multi-dimensional kernel convolution, the input function $a(x)$ is projected into a high-dimensional space  $d_v$ using a local transformation $\mathrm{P}: a(x) \mapsto v_0(x)$. 
This local transformation $\mathrm{P}$ can be implemented as a shallow, fully connected neural network (FNN). Within the lifted space, $l$  iterations  equivalent of the form  similar to \autoref{eq:hammerstein} are executed, $G: \mathbb{R}^{d_v} \mapsto \mathbb{R}^{d_v}$ such that  $v_{j+1} = G(v_j)$. Once $l$-iterations are completed, the output undergoes another local transformation ${\rm{Q}}: v_{l}(x) \mapsto u(x)$ to obtain the final solution $u(x) \in \mathbb{R}^{d_u}$. 
The step-wise iteration employing the composition operator, $G(\cdot)$ can be mathematically expressed as follows:
\begin{equation}\label{eq:iteration}
    G(v_{j+1})(x):= \varphi \left( \left(K(a; \phi) * v_{j}\right)(x) + W v_{j}(x) \right), \quad x \in D, \quad j \in [1,l], 
\end{equation} 
where $\varphi(\cdot) \in \mathbb{R}$ is a non-linear activation operator, $\bm \phi \subset \bm{\theta}_{NN}$ represents the trainable kernel parameters, $W: \mathbb{R}^{d_v} \to \mathbb{R}^{d_v}$ is a linear transformation, and $K$ denotes the nonlinear integral operator given by:
\begin{equation}\label{eq:integral}
\left(K(a; \bm \phi) * v_{j}\right)(x) := \int_{D} k \left(a(x), x, \xi; \phi \right) v_{j}(\xi) \mathrm{d}\xi, \quad x \in D, \quad j \in [1,l].
\end{equation} 
In the above expression, $k \left(a(x), x, \xi; \phi \right)$ represents the kernel of the nonlinear integral equation in \autoref{eq:hammerstein}. These kernel parameters, denoted as $ \bm \phi$, are learnable and obtained through end-to-end training of the neural operator. In WNO, it is proposed to replace the integral operator \autoref{eq:integral} by a convolution operation defined in the wavelet transformed space. To that end,  a wavelet transform is performed on the lifted input $v_{j}(x)$. The forward wavelet transform,$(\mathcal{W}(\cdot)$ and inverse wavelet transform $\mathcal{W}(\cdot)$, are defined \cite{daubechies1992ten} as follows:
\begin{equation}\label{eq:wavelet}
    \begin{aligned}
        (\mathcal{W} v_{j})(s, \tau) & = \int_{D} \Gamma (x) \frac{1}{|s|^{1 / 2}} \psi\left(\frac{x-\tau}{s}\right) dx, \\
        (\mathcal{W}^{-1} (v_{j})_w)(x) & = \frac{1}{C_{\psi}} \int_{0}^{\infty} \int_{D} (v_{j})_{w}(s, \tau) \frac{1}{|s|^{1 / 2}} \tilde{\psi}\left(\frac{x-\tau}{s}\right) d\tau \frac{ds}{s^{2}},
    \end{aligned}
\end{equation}
Here, $\psi(x)$ represents the orthonormal mother wavelet, and $s$ and $\tau$ are the scaling and translational parameters used in the wavelet decomposition. While $(v_{j})_{w}$ denotes the decomposed wavelet coefficients of $v_{j}(x)$, $\psi(\cdot)$ refers to the mother wavelet that is scaled and shifted. $ C_{\psi}$ is the admissible constant such that $0 < C_{\psi} < \infty$ \cite{daubechies1992ten}.
The integral kernel \(k_{\phi}\), defined in the wavelet space, is denoted as \(R_{\phi} = \mathcal{W}(k_{\phi})\). Now, using the convolution theorem for parameterized kernels, the kernel integration in \autoref{eq:integral} carried out over the wavelet domain is expressed as:
\begin{equation}\label{eq:conv_final} 
\left(K(\phi) * v_{j}\right)(x) = \mathcal{W}^{-1}\left(R_{\phi} \cdot \mathcal{W}(v_{j})\right)(x), \quad x \in D. 
\end{equation}
Executing complex wavelet decomposition in \autoref{eq:wavelet} is computationally expensive as the scale and translation parameters $s$ and $\tau$ are theoretically infinite-dimensional. However, this operation can be streamlined by utilising the slim dual-tree complex wavelet transform (DTCWT) toolbox \cite{cotter2020uses}, originally proposed in \cite{selesnick2005dual}. Slim DTCWT provides two sets (real and imaginary) of six wavelet coefficients at each level of decomposition, which approximates the coefficients of 15$^{\circ}$, 45$^{\circ}$, 75$^{\circ}$, 105$^{\circ}$, 135$^{\circ}$, and 165$^{\circ}$ wavelets. Here, it should be noted that the coefficients are halved at each decomposition level due to its conjugate symmetry. The most relevant information of the input is obtained at the highest level of the DTCWT, where the coefficients corresponding to the lower frequencies are preserved. Thus, parameterization of the kernel $R_{\phi}(\cdot)$ is done at the same level. The overall kernel convolution in the wavelet domain,  $(R_{\phi} \cdot \mathcal{W}( v_{j}))(x)$ can be expressed as follows,
\begin{equation}\label{eq:convolution_final}
    \begin{aligned}
        \left(R \cdot \mathcal{W} (v_{j}; \ell)\right)_{t_1,t_2} = \sum_{t_{3}=1}^{d_{v}} R_{t_1, t_2, t_3} \mathcal{W} (v_{j}; \ell)_{t_1, t_3}; && I_1 \in [1, {d_l}], && I_{2}, I_{3} \in d_{v}.
    \end{aligned}
\end{equation}
Here, ${d_l}$ denotes the dimension of each wavelet coefficient at the last level of DTCWT, where the wavelet transform is applied on the uplifted input $v_{j}(x)$. To obtain the decomposed output $\mathcal{W}\left(v_{j}; \ell \right)$ having dimension ${{d_l} \times d_{v}}$ we construct the weight tensor $R_{\phi}(\ell)$ with the dimension ${d_{l} \times d_{v} \times d_{v}}$. The DTCWT contains 12 wavelet coefficients, including real and imaginary coefficients of 15$^{\circ}$, 45$^{\circ}$, 75$^{\circ}$, 105$^{\circ}$, 135$^{\circ}$, and 165$^{\circ}$ wavelets. Thus, in order to learn the parametric space, we need 12 weight tensors and corresponding 12 convolutions (\autoref{eq:convolution_final}). 

\begin{figure}[!ht]
    \centering
    \includegraphics[width=\textwidth]{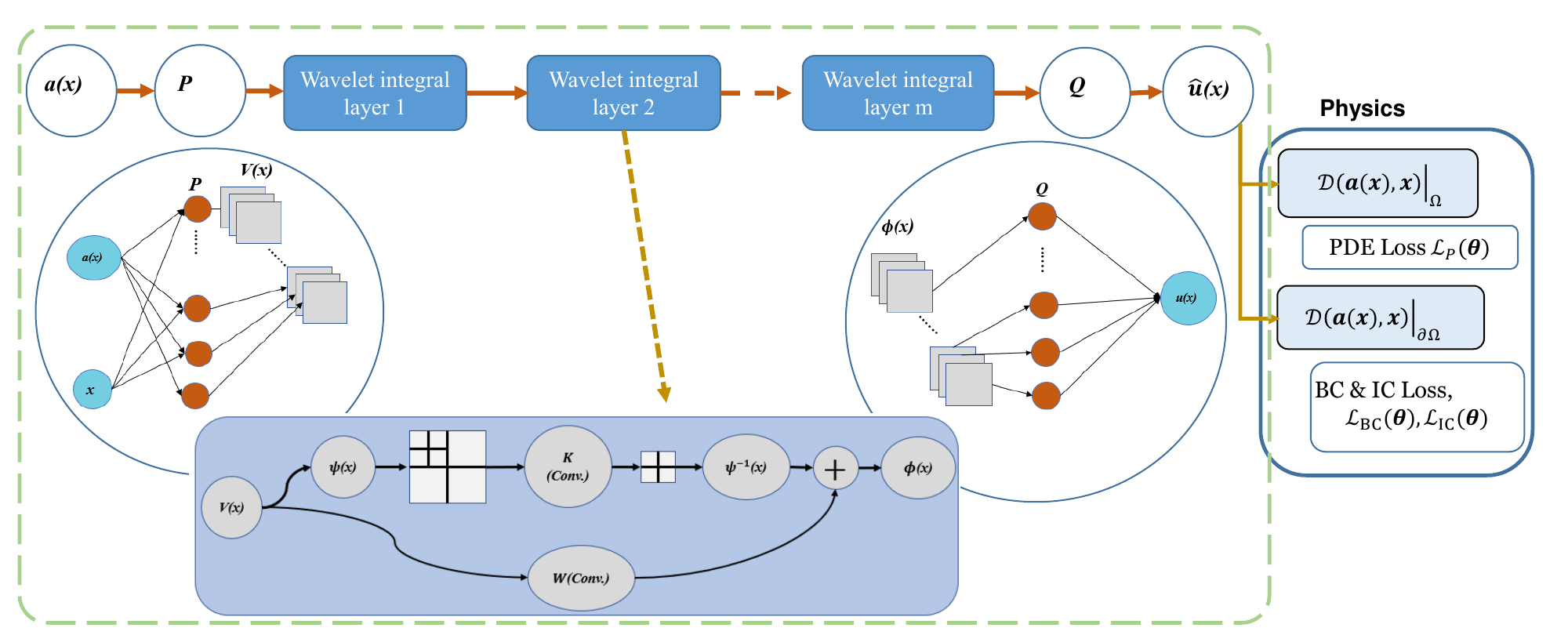}
    \caption{\textbf{The Proposed Physics-informed Operator (PIO) for reliability analysis}. Here, we propose physics-informed WNO (PIWNO) as the PIO, where inputs are initially lifted to a high-dimensional latent space, where they undergo iterative processes. These iterations are represented using wavelet kernel integration blocks, which consist of a kernel integration network that learns the integration kernel and a linear transformation network. The latent inputs are transformed into the space-frequency localised domain using wavelets. The outputs of the integration and the integral constants are then combined, and a nonlinear activation is applied. Output from the ultimate integral layer is down-lifted to obtain the final output, which is obtained as the solution for the underlying PDE. The solutions are constrained to satisfy the given PDEs, boundary conditions (BC), and initial conditions (IC). To enforce the PDE constraints, spatial derivatives are computed using a stochastic projection-based gradient estimation scheme.}
    \label{fig:piwno}
\end{figure}
To summarize, the Wavelet Neural Operator (WNO) is mainly defined by four hyperparameters: (i) uplifting dimensions, (ii) vanishing moments in the wavelets, (iii) decomposition level, and (iv) the number of wavelet blocks. The uplifting dimension can be seen as increasing the channel dimension in a convolutional neural network (CNNs), which facilitates learning. On the other hand, vanishing moments define the smoothness of the wavelet, such that the wavelet with a lower order moment is suitable for capturing images with higher spatial variations, and the wavelet with a higher order moment captures smoothly varying images effectively. Kernel size can be decided based on the number of levels of decomposition, as at each level of decomposition, the input is sub-sampled by a factor of two. Thus, the optimum number of kernel parameters is chosen based on the level of decomposition. Most importantly, the number of wavelet blocks selected is based on the complexity of the underlying operator to be learned. All these parameters are to be decided based on the validation set or by employing a neural architecture search algorithm such as the one presented in \cite{soin2024generative}.

\subsection{Physics informed WNO}
\begin{figure}
    \centering
    \includegraphics{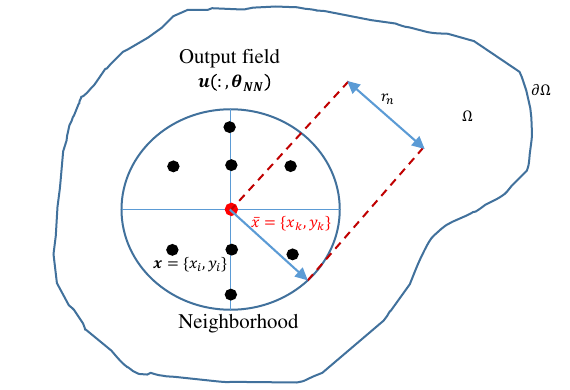}
    \caption{A diagrammatic representation grid point (red dot) and neighbourhood region utilized compute stochastic projection based gradient, where the black dots denote the neighbourhood collocation points}
    \label{fig:sp_projection}
\end{figure}

In a data-driven setting, an operator-learning framework $\hat{\mathcal{M}}$ can be trained using $N$ pairs $\{a_j, u_j\}_{j=1}^{N}$ to approximate the operator $\mathcal{M}$. However, in scenarios where observed output data is unavailable, or data generation requires expensive simulations, we rely on the underlying physics to train WNO using \textit{only} the inputs $\{a_j\}_{j=1}^{N}$.
To implement a physics-informed operator, we revisit \autoref{operator1}, where differential operator $\mathcal{N}$, includes all the derivatives with respect to space and time (i.e. $\partial_{\mathbf{t}}, \partial_{\mathbf{t}}^2, \ldots,  \partial_t^n, \partial_x, \ldots, \partial_x^n $). The residual form of differential equation \autoref{operator1} can be rewritten as: 
\begin{equation}\label{diff: eq}
\mathcal{N}\left(\bm{a},\bm{u},\bm{x}, t\right) - \bm{s}(\bm x,t) = 0. 
\end{equation}
Training physics-informed operators requires incorporating residual loss calculated from the governing PDE. For computing the loss, the output prediction from the WNO, i.e. $\hat{\bm u} = \hat{\mathcal{M}}(\bm{a}(\bm x), \bm{\theta}_{NN})$, is plugged into the \autoref{diff: eq}. Subsequently, the total physics loss is evaluated by summing the PDE loss, boundary loss, and initial condition loss with each component appropriately scaled.
The expression of physics loss (also known as total residual loss), $\mathcal{L}_{\text{Physics}}$ is given by the following: 
\begin{equation}\label{Lpde} 
\mathcal{L}_{\text {Physics }}\left(\bm{\theta}_{NN}\right) = \underbrace{\left\|\mathcal{N}(\bm{a}(\bm x,t
), \hat{\bm{u}})-\bm{s}(\bm x,t)\right\|^2}_{\text{PDE Loss}} + {\alpha_1} \underbrace{ \left\|{\hat{\bm{u}}(\bm{x},t) \mid}-\bm{g}\right\|^2}_{\text{Boundary Loss}} + {\alpha_2} \underbrace{ \left\|{\hat{\bm{u}}(\bm{x},t=0) \mid}-\bm{\hat{\bm{u}}_{0}}\right\|^2}_{\text{Initial Condition Loss}}
\end{equation}
For $N$ training input samples, with $n_d$ space discretizations, $n_b$ boundary points,  and $n_{ic}$ initial condition points \autoref{Lpde} can be rewritten as: 
\begin{equation}\label{sp_projection}
\begin{aligned}
\mathcal{L}_{\text{Physics}}\left(\bm{\theta}_{NN}\right) = \frac{1}{N} \sum_{j=1}^N \sum_{i=1}^{n_d}\left|\mathcal{N}(a_j, \hat{u}_j)(x_i,t)-s(x_i,t)\right|^2 + & {\alpha_1} \sum_{j=1}^N \sum_{i=1}^{n_b}\left|\hat{u}_j(x_i,t) - g(x_i,t)\right|^2+ \\
&{\alpha_2} \sum_{j=1}^N \sum_{i=1}^{n_{ic}}\left| \hat{u}_j(x_i,t=0) - u_0(x_i)\right|^2 
\end{aligned}
\end{equation}
The constants $\alpha_1$ and $\alpha_2$ in the above equation represent scaling weights given to boundary and initial conditions losses. Finally, minimizing the given loss function yields the optimal network parameters:
\begin{equation}
    \bm{\theta}^{*}_{NN} = \underset{\bm{\theta}_{NN}} {\text{argmin}}\; \mathcal{L}_{\text {Physics}}\left(\hat{\mathcal{M}}(\bm{a}(\bm x);\bm{\theta}_{NN})\right).
\end{equation}
It is to be noted here that since the operator is trained with loss functions computed directly from the governing equations $\mathcal{L}_{Physics}$, one needs to evaluate the derivatives operator output involved in the governing physics. However, in practice, it is non-trivial to obtain the derivatives of the output field,  especially when the architecture of  $\hat{\mathcal{M}}({\cdots;\bm{\theta}}_{NN})$ is devised of the convolution layers. Thus, we conveniently employ stochastic projections to compute the gradients involved in the training loss. To compute gradients at a given point $\bar{x} = (x_k, y_k)$ within the domain using stochastic projection \cite{navaneeth2023stochastic}, we consider a neighbourhood within a specified radius $r_n$ (see \autoref{fig:sp_projection}). Within this neighbourhood, $N_t$ collocation points are chosen. The gradient of the output $\hat {\bm u}$ with respect to the input variable at $\bar{x}$ is calculated using the expression given by:
\begin{equation}\label{spgf_net}
    \mathbf{\hat G}(\bm x={\bm{\bar{x}}})  = \frac{\partial \hat{\bm{u}}({\bm {\bar x}})}{\partial {\bm {x}}}= \frac{\frac{1}{N_b}{\sum}_{i=1}^{N_{b}}{(\hat{\bm{u}}(\bm {x}_i)- \hat{\bm{u}}({\bm \bar{x}}))(\bm {x}_i-\bm {\bar {x}})^{T}}}{\frac{1}{N_b}{\sum}_{i=1}^{N_{b}}{(\bm {x}_i-\bm {\bar{x}})(\bm {x}_i-\bm {\bar{x}})^{T}}}
\end{equation}
where $\bm{x}_i= \{x_i,y_i\}$ is considered to be a generic neighborhood point and $N_b$ represents the number of neighborhood points.
The step-by-step training procedure for the Physics-Informed Operator (PIO) is detailed in Algorithm \ref{alg:Training PIO}. Additionally, the implementation steps for computing reliability using the trained PIO are outlined in Algorithm \ref{alg:Reliabiilty}.

\begin{algorithm}[ht!]
    \caption{Training algorithm for physics informed operator}\label{alg:Training PIO}
    \textbf{Requirements:} Boundary conditions, initial conditions, source functions and PDE describing the physics constraint.
    \begin{algorithmic}[1]
    \State {\textbf{Initialize:} Network parameters, $\bm{\theta}_{NN}$ of the WNO, $\hat{\mathcal{M}}$.
    
    \State Pass the input $\bm{a(x)}$ to WNO and the grid points over the domain $\{x^{i}_{f}, y^{i}_{f} \} \in \Omega$.
    
    \State Obtain the output prediction of WNO ($\hat {\mathcal U}$) corresponding coordinates of the grid points over the domain $\{x^{i}_{f}, y^{i}_{f} \} \in \Omega$.
    
    \State Collect the output prediction of WNO ($\hat{\mathcal U}$) and corresponding coordinates of the grid points at the boundary of the domain $\{x^{i}_{b}, y^{i}_{b} \} \in \partial \Omega$.
    
    \State For the given resolution of the field, define the neighbourhood of each grid point based on the radius $r_n$.
    
    \State Obtain the gradients of the field variable at all the collocation points using Eq.\ref{sp_projection} and store the gradients. 
    
    \State Compute the PDE loss $\mathcal{L}_{PDE}$ using the gradeint obtained from the previous step.
    
    \State Compute the boundary loss $\mathcal{L}_{BC}$, sum the all losses to get the total loss $\mathcal{L}_{Physics}$.
    
    \While {{$\mathcal{L} > \epsilon$}}
    \State {Train the network:} $\bm {\theta}_{NN}\leftarrow \bm{\theta}_{NN}-\delta \nabla_{{\bm{\theta}}_{NN}}{\mathcal L}(\bm \theta_{NN})$.
    \State {epoch= epoch $+$ 1.}
    \EndWhile
    \State {Return the optimum parameters for the PIWNO ($\bm{\theta}^{*}_{NN}$).}}
    \end{algorithmic}
    {\textbf{Output:} Trained PIO}
\end{algorithm}

\begin{algorithm}[ht!]
    \caption{Physics Informed Operator for reliability analysis.}\label{alg:Reliabiilty}
    \textbf{Requirements:} Trained PIO, limit-states for the given problem, and samples of the input conditions.
    \begin{algorithmic}[1]
    \State {Load the trained PIO ($\hat{\mathcal{M}}(\cdots;\bm{\theta}^{*}_{NN})$).}
    \State {Obtain predictions/solutions.}
    \For{$i=1,\ldots,N_s$}
    \State {Draw $i-$th sample from the distribution of the input functions}
    \State {Evaluate the limit state function.}
    \EndFor
    \State {Obtain Probability density function/histogram of the quantity of interest, Failure probability.}
    \end{algorithmic}
    {\textbf{Output:} Probability density function, Failure probability and Reliability Index}
\end{algorithm}

\section{Numerical examples}\label{sec: Numerical example}
We validate the efficacy of the Physics-Informed Operator (PIO) through detailed illustration of four numerical examples in this section. As was mentioned, the overarching goal is to quantify the reliability of a parametric system with varying inputs. Here, the varying input functions for the chosen numerical examples include initial conditions, source functions and input parametric fields. It is emphasised that only input samples are needed to train the physics-informed operator, where the numerical examples include the stochastic field with an intrinsic dimensionality of up to 358.
The examples selected involve both time-dependent and time-independent problems. The failure probability ($P_f$) is obtained based on the limit state function defined for individual problems, and the corresponding reliability index ($\beta$) is given by:
\begin{equation}
    \beta=\Phi^{-1}(1-P_f)
\end{equation}
Here, we note that for the time-dependent problems, the Probability Density Function (PDF) is obtained based on the first passage failure time. The architecture of the PIO framework is comprised of three to five layers, depending on the specific example under consideration. ADAM optimizer, with a weight decay of $10^{-6}$ and an initial learning rate of 0.001, is employed to optimize the WNO parameters. The batch size varies from 10 to 25, while the overall number of epochs is fixed to 300.

\subsection{Diffusion-reaction system}
The first example we consider is the non-linear diffusion-reaction system. Diffusion-reaction process is a well-studied phenomenon in chemical systems that involves the movement of chemicals through space.
One such practical scenario is the spread of the pollutants introduced into the water, which undergoes diffusion due to the natural movement of water molecules and simultaneously undergoes chemical reactions, such as decay or neutralization. The reaction-diffusion equation describes the mathematical model for the physical phenomena. The limit state here can be defined in terms of the maximum concentration of the chemical substance/pollutant at a given spatial point. More precisely, the limit state function defines the failure that occurs if the maximum concentration ($u$) at a given spatial location, defined in terms of the spatial grid, such that $x_{grid} = x_{sp}$ crosses threshold $e_h$. Mathematically, it is defined as:
\begin{equation}
    \mathcal{J}(x)=e_h-{|u(x_{sp},t)|}_{max},
\end{equation}
where $u(x,t)$ is obtained by solving the following equation,
\begin{equation}\label{eq:eg1}
    \begin{aligned}
    \partial_t u -{B}\partial_{x x} u -{k}{u}^2&= f(x), & & x \in(0,1), t \in(0,1] \\
    u(x=0, t) &=u(x=1, t)=0, & & x \in(0,1), t \in(0,1] \\
    u(x, 0) & = 0, & & x \in(0,1).
    \end{aligned}
\end{equation}
We have considered $B=0.01$ and $k=0.01$. To quantify the reliability of the system, we need to compute the output response. For that, we seek an operator mapping source function $f(x)$ to the corresponding output, $u(x,t)$, the spatiotemporal solution of the PDE. For reliability analysis, we have considered the forcing function $f(x)$ as random and modelled it as follows,
\begin{equation}\label{eq:input1}
 f(x) = {n} \sin(\pi x) + (1 - {n}) \cos(\pi x) + {p} \sin(2 \pi x) + (1 - {p}) \cos(2 \pi x) + {w} \sin(3 \pi x) + (1 - {w}) \cos(3 \pi x),
\end{equation}
where the parameters $n$,$p$ and $w$ randomly chosen from uniform distributions such that
$n \sim \operatorname{Unif}(0, 1)$, $p \sim \operatorname{Unif}(0, 1)$ and $w \sim \operatorname{Unif}(0, 1)$. This indicates that this is a low-diemnsional problem with intrinsic dimensionality of 3. For training the physics-informed operator, we generated 600 random realizations on a fixed resolution $81 \times 81$ of the random field $f(x)$ by using Eq. \eqref{eq:input1}. Note that as the training is based on the governing physics in Eq. \eqref{eq:eg1}, no output data is required. 

\begin{figure}[!h]
    \centering{
    \includegraphics[width=1\textwidth]{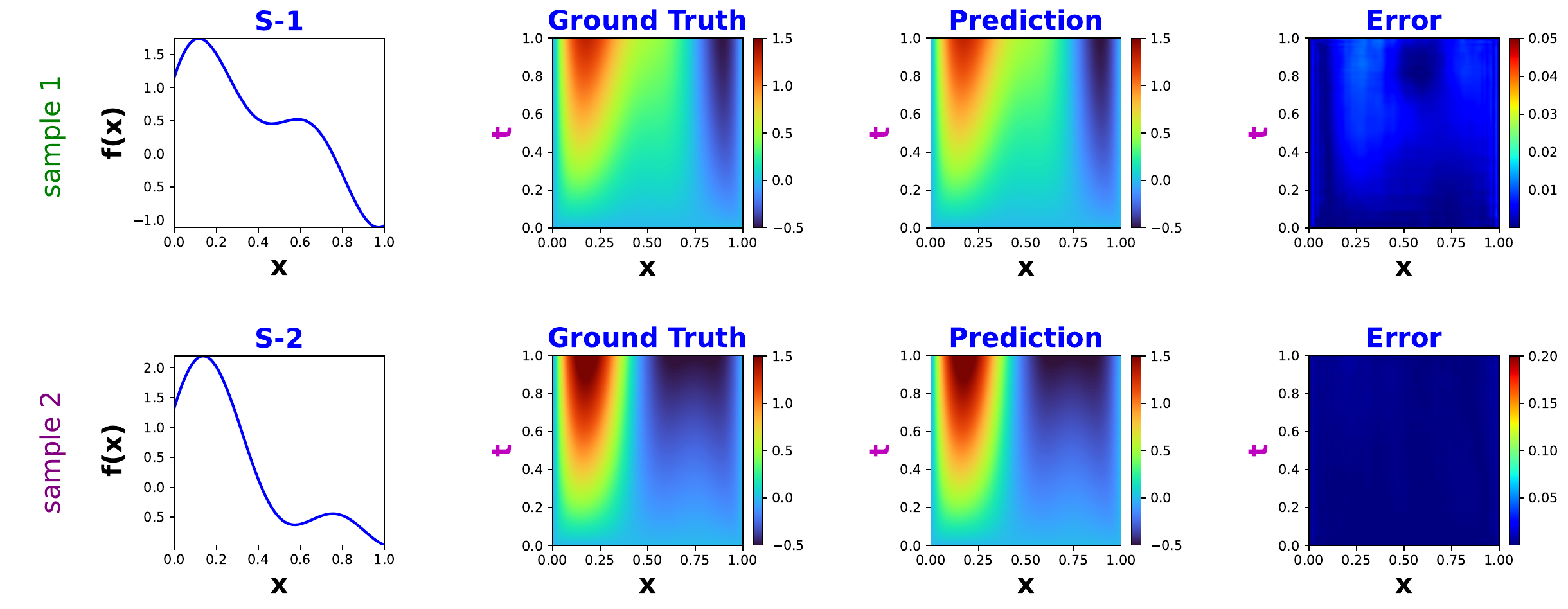}}
    \caption{The results for the Diffusion-reaction system comprised of source functions, ground truth solutions, predictions, and error plots demonstrated for using 2 different unseen sample instances. The PIO effectively maps the initial condition to the corresponding solution $u(x,t)$ over the domain, with a spatiotemporal resolution of $81 \times 81$.}
    \label{fig:reactdiffphy}
\end{figure}

\begin{figure}[!h]
    \centering{
    \includegraphics[width=0.8\textwidth]{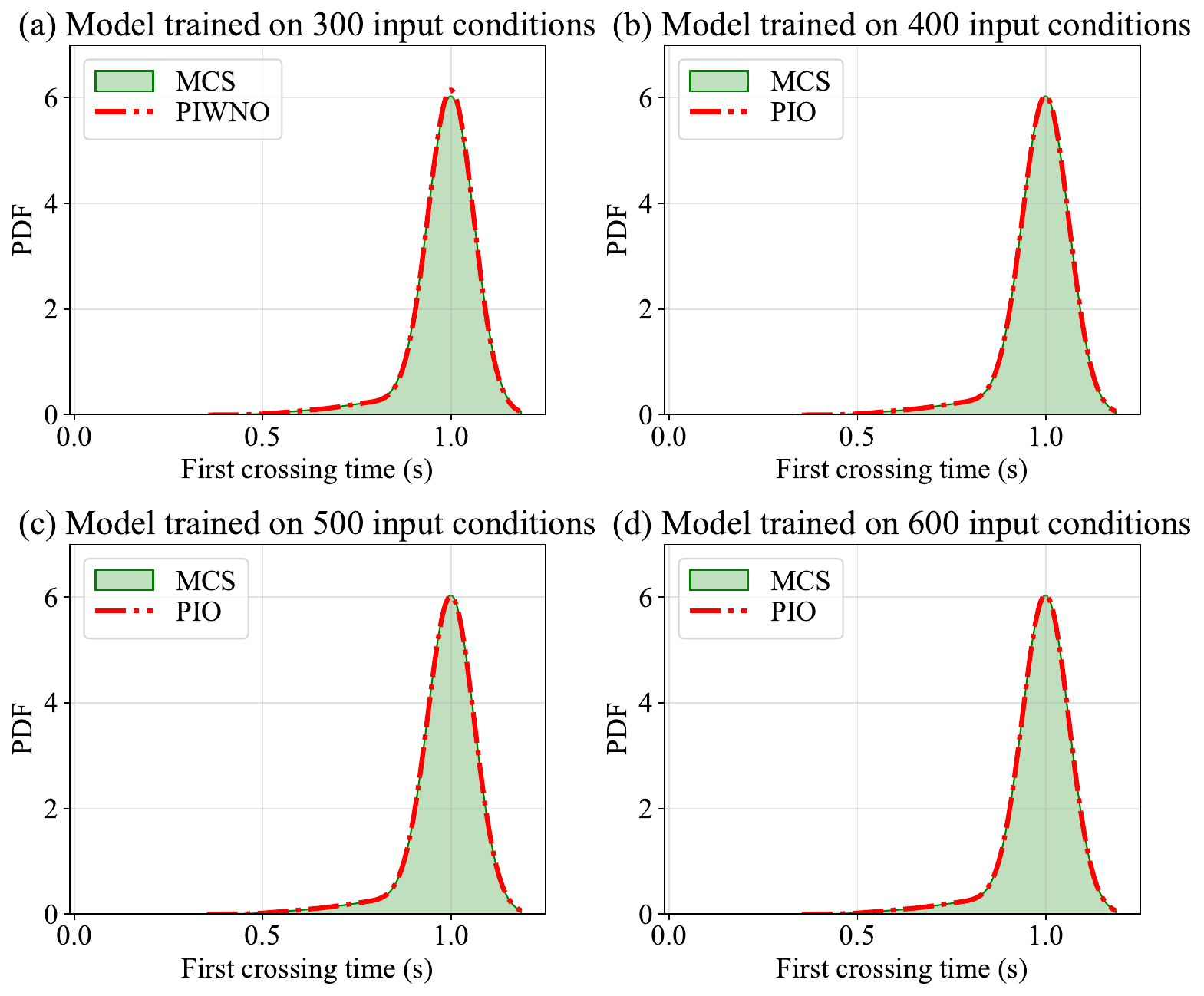}}
    \caption{PDF plots of failure time obtained by MCS and PIO trained with varying number of source functions for the diffusion-reaction system}
    \label{fig:reactdiffpdf}
\end{figure}

We first validate the predicted response obtained using PIO with those obtained using numerically solving the equation. The contours presented in \autoref{fig:reactdiffphy} indicate that the proposed approach matches well with the ground truth. Reliability analysis is carried out based on the limit state function aforementioned, where the grid is chosen to be $x_{sp} = 41$ and the threshold to be $e_h =0.85$. Subsequently, the failure probability ($p_f$) is defined as the ratio of the number of samples that fail to satisfy the limiting state to the total number of samples. Probability density function corresponding to the first passage failure is illustrated in  \autoref{fig:reactdiffpdf}. Results corresponding to different numbers of training inputs are reported. For benchmarking, results using vanilla Monte Carlo Simulation (MCS) are also generated. It is observed that PIO yields highly accurate results for all the cases. Estimated failure probability and the reliability index are provided in the \autoref{table1}, which indicates the efficacy of the PIO in accurately computing the failure probability and reliability. For a comprehensive analysis, the variation of the failure probability with the threshold is also computed using MCS, first-order reliability method (FORM), second-order reliability method (SORM), data-driven WNO (with varying training samples), and PIO. It is observed that PIO and SORM yield the best result followed by the data-driven WNO. FORM, on the other hand, fails to capture the probability of failure accurately. This can perhaps be attributed to the nonlinearity present in Eq. \eqref{eq:eg1}.

\begin{figure}[!h]
    \centering{
    \includegraphics[width=0.8 \textwidth]{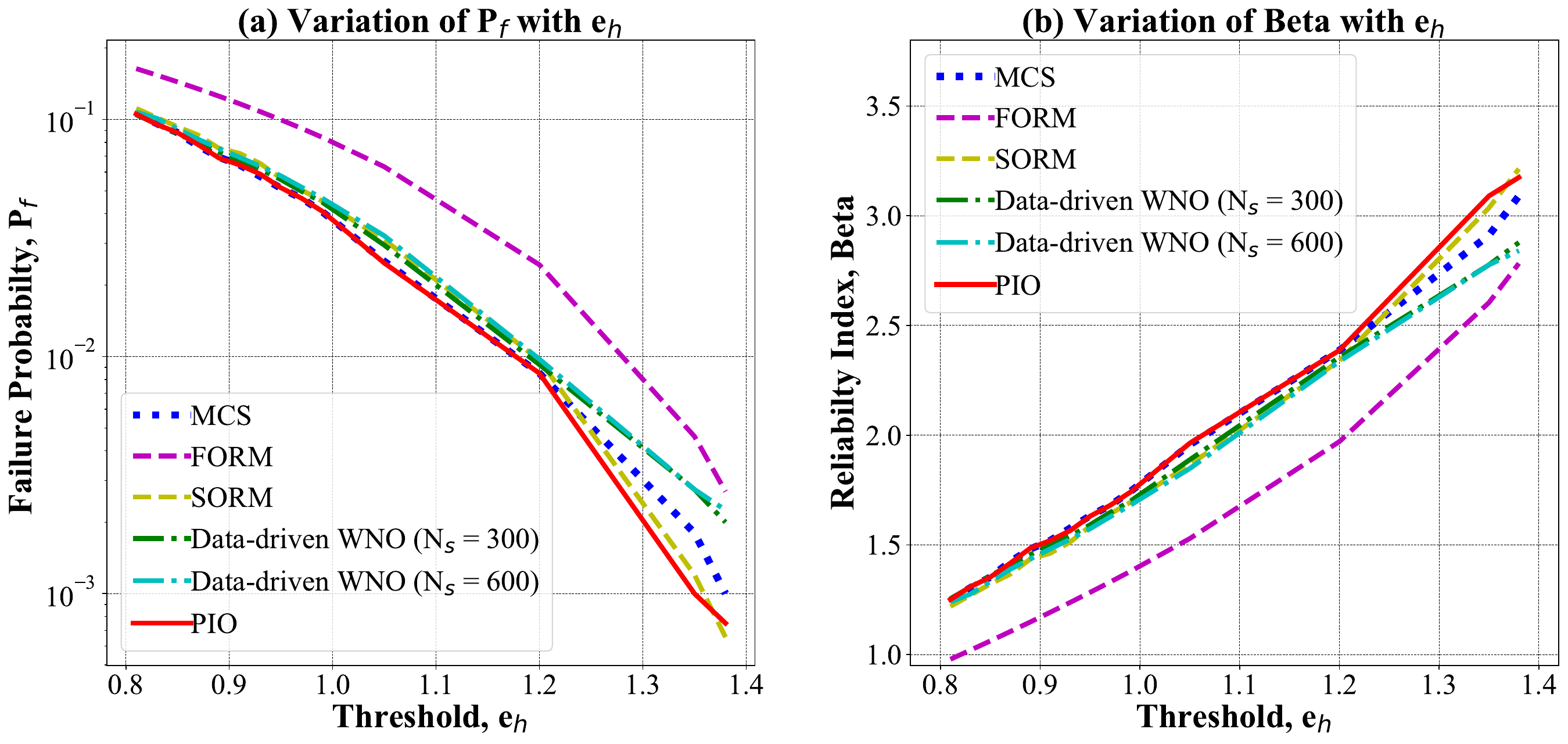}}
    \caption{Prediction results of the diffusion-reaction system for the failure probability ($P_f$) and the reliability index ($\beta$) with unceasing limit state threshold, obtained by PIO in comparison with the results of MCS, FORM, SORM and data-driven WNOs (trained with a number of samples, $N_s = 300$ and $N_s = 600$) }
    \label{fig:reactdiff_vare_h}
\end{figure}
\begin{table}[ht!]
    \centering
    \caption{Results of first passage failure probability obtained using   
    proposed framework for the diffusion-reaction system with varying trained source functions.}
    \label{table1}
\begin{tabular}{lccccc} 
\hline
\textbf{No. of training conditions}& \textbf{300} & \textbf{400} & \textbf{500}&\textbf{600}&\textbf{Actual}\\
\hline
Failure probability ($P_f$) &  0.0855  &  0.0887 &   $0.0885$ & 0.0880 & 0.0880\\ \hline
Reliability index ($\beta$) &  1.369 &  1.348  &  1.350 & 1.353 & 1.353\\ 
\hline 
\end{tabular}
\end{table}

\subsection{Impulse transmission in nerve}
\begin{figure}[!h]
    \centering{
    \includegraphics[width=0.4 \textwidth]{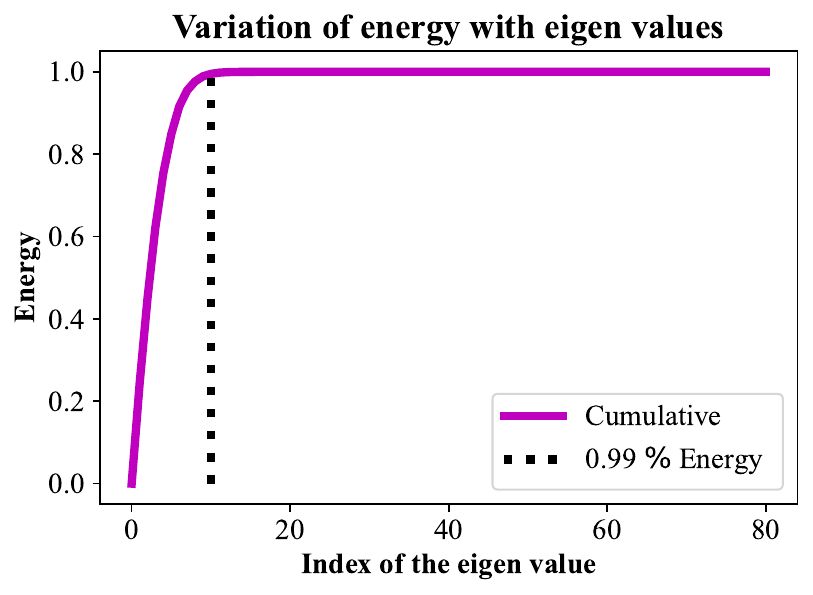}}
    \caption{Intrinsic dimensionality of the input function obtained for example of nerve impulse transmission, where the 99$\%$ of the energy is contained in first 10 eigenvalues.}
    \label{fig:Nagumo_Intrinsic}
\end{figure}

As the next example, we consider a problem involving impulse transmission in the nervous system. This problem involves propagation of electrical signals through neurons and is crucial for understanding various functions, including muscle contraction, sensory perception, and cognitive activities. The driving reason for transmitting impulses in human nerve fibres is the dynamic electric potential changes across nerve cells. The mechanism can be effectively modelled using the Nagumo equation \cite{mckean1970nagumo}. During the transmission, exceeding this action potential above the typical range can affect the repolarization of the axon membrane, disrupt normal nerve function and lead to failure of impulse transmission. Therefore, the limit state function is defined in terms of the maximum electrical potential field generated over the space and given time. The failure of the system is considered when the maximum absolute value of response ($u$) at the spatial point, grid $x_{grid} = x_{sp}$ crosses the threshold $e_h$. Thus, the limit-state function is represented as:
\begin{equation}
    \mathcal{J}(x)=e_h-{|u(x_{sp},t)|}_{max},
\end{equation}
where $u(x,t)$ is obtained by solving the Nagumo equation \cite{laing2009stochastic,lord2014introduction},
\begin{equation}
    \begin{aligned}
    \partial_t u - \varepsilon \partial_{x x} u &= u(1-u)(u-\alpha), & & x \in(0,1), t \in(0,1] \\
    u(x=0, t) &=u(x=1, t)=0, & & x \in(0,1), t \in(0,1] \\
    u(x, 0) & =u_0(x), & & x \in(0,1) ,
    \end{aligned}
\end{equation}
where the parameter $\varepsilon > 0$ determines the rate of diﬀusion and $\alpha \in \mathbb{R}$ decides the speed of a wave travelling down the length of the axon. We have considered $\varepsilon=1$ and $\alpha=-1/2$. The objective here is to quantify the reliability of this system due to uncertainty in the initial condition, which in turn is modeled as a Gaussian random field with a specified kernel given below:
\begin{equation}
\mathcal{K}(x, y)=\sigma^{2} {exp}\left(\frac{{-}(\bm x-\bm{x}^{'})^{2}}{2{l}^{2}}\right).
\end{equation}
We have considered $\sigma= 0.1;$  and $l = 0.1$. We employ Karhunen-Lo{\`e}ve expansion (KLE) to compute the intrinsic dimensionality of the problem. The energy plot in \autoref{fig:Nagumo_Intrinsic} indicates that the first ten leading Eigenvalues capture 99\% of the energy, and hence, it is also a relatively low-dimensional problem. PIO is trained to map the initial conditions $u_0(x)$ to spatiotemporal solutions $u(x, t)$, i.e., $\mathcal{M}: u_0(x) \mapsto u(x, t)$. The network is trained using 800 randomly generated initial conditions on $65\times65$ grid.

\begin{figure}[!h]
    \centering{
    \includegraphics[width=1\textwidth]{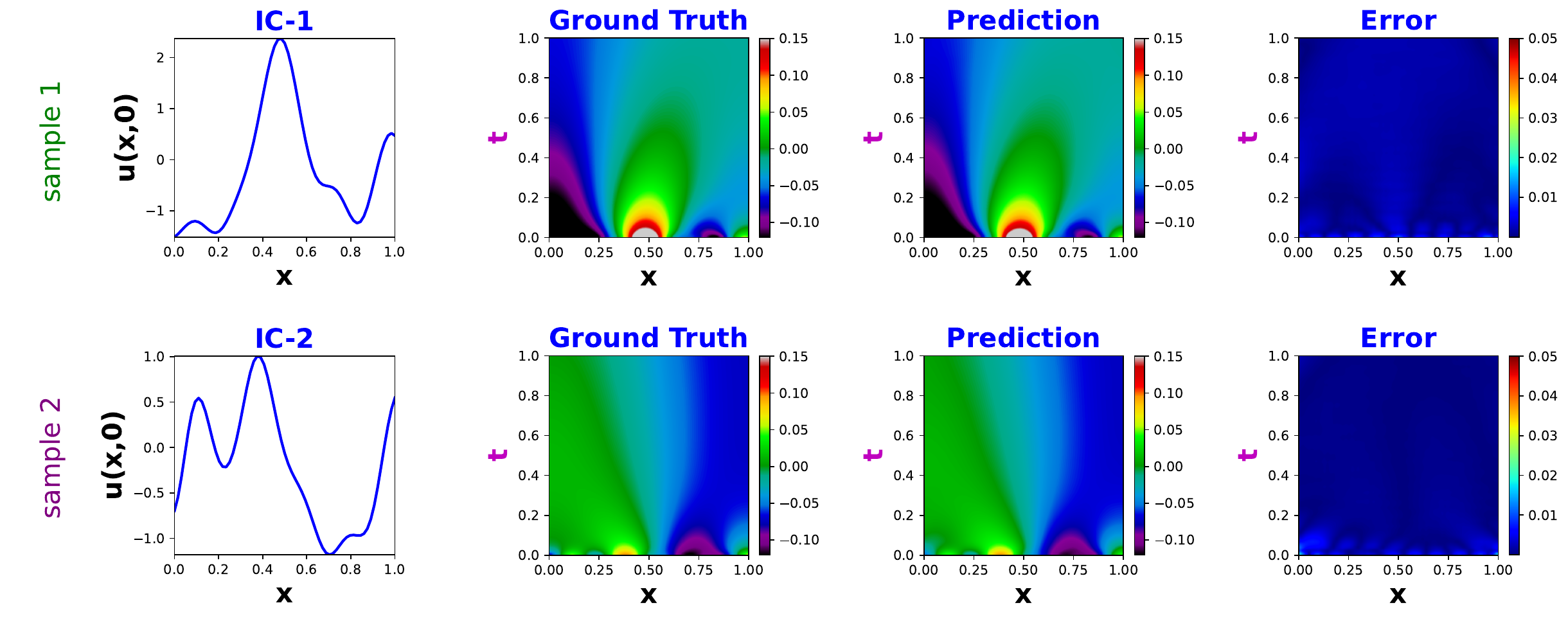}}
    \caption{The results for the nerve impulse transmission example comprised of initial conditions, ground truth solutions, predictions, and error plots demonstrated for using  2 different unseen sample instances. The PIO effectively maps the initial condition to the corresponding solution $u(x,t)$ over the domain, with a spatiotemporal resolution of $65 \times 65$.}
    \label{fig:nagumophy}
\end{figure}

\begin{figure}[!h]
    \centering{
    \includegraphics[width=0.8\textwidth]{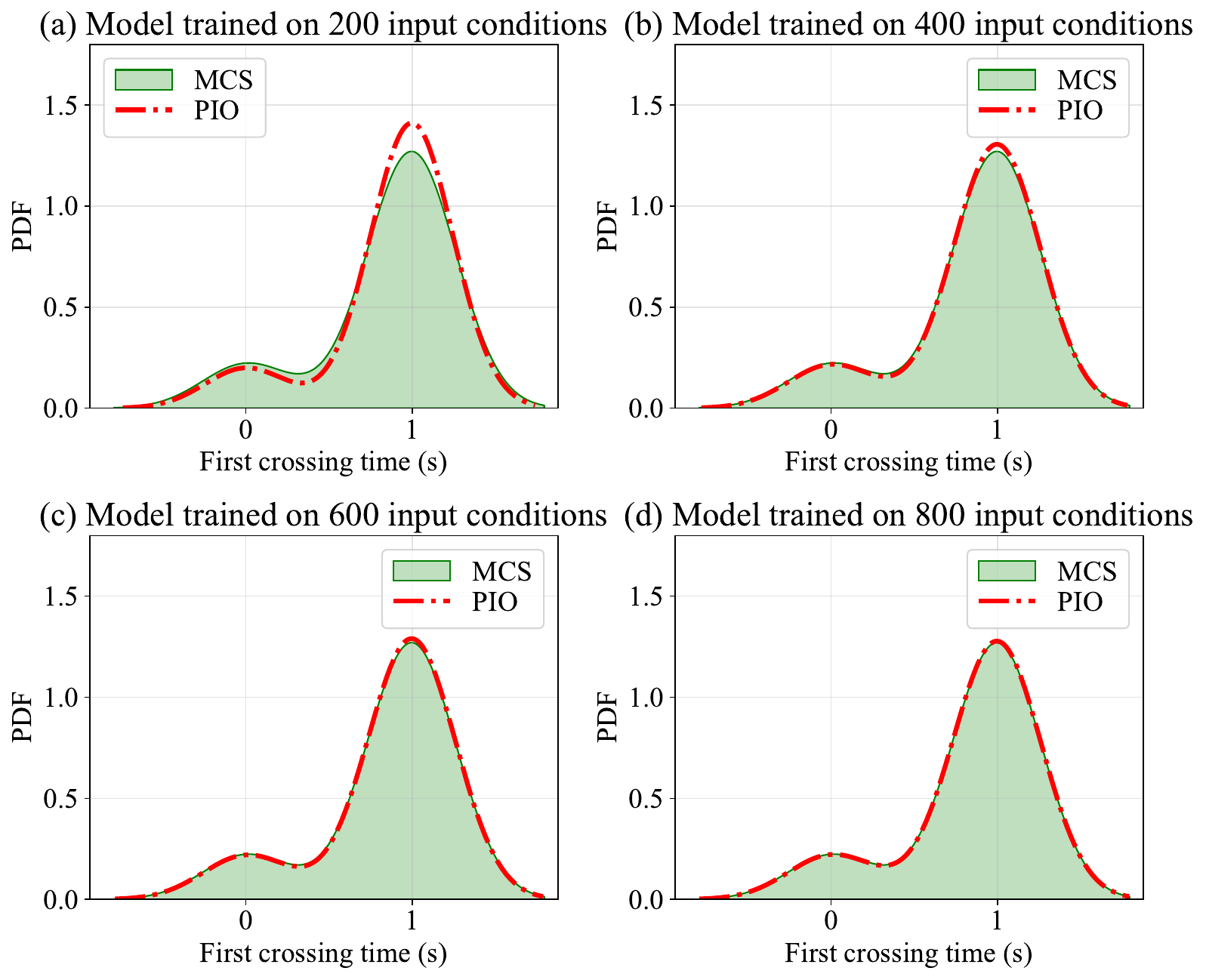}}
    \caption{PDF plots of failure time obtained by MCS and PIO trained with varying number of initial conditions for the nerve impulse transmission example}
    \label{fig:nagumopdf}
\end{figure}

\begin{figure}[!h]
    \centering{
    \includegraphics[width=0.8 \textwidth]{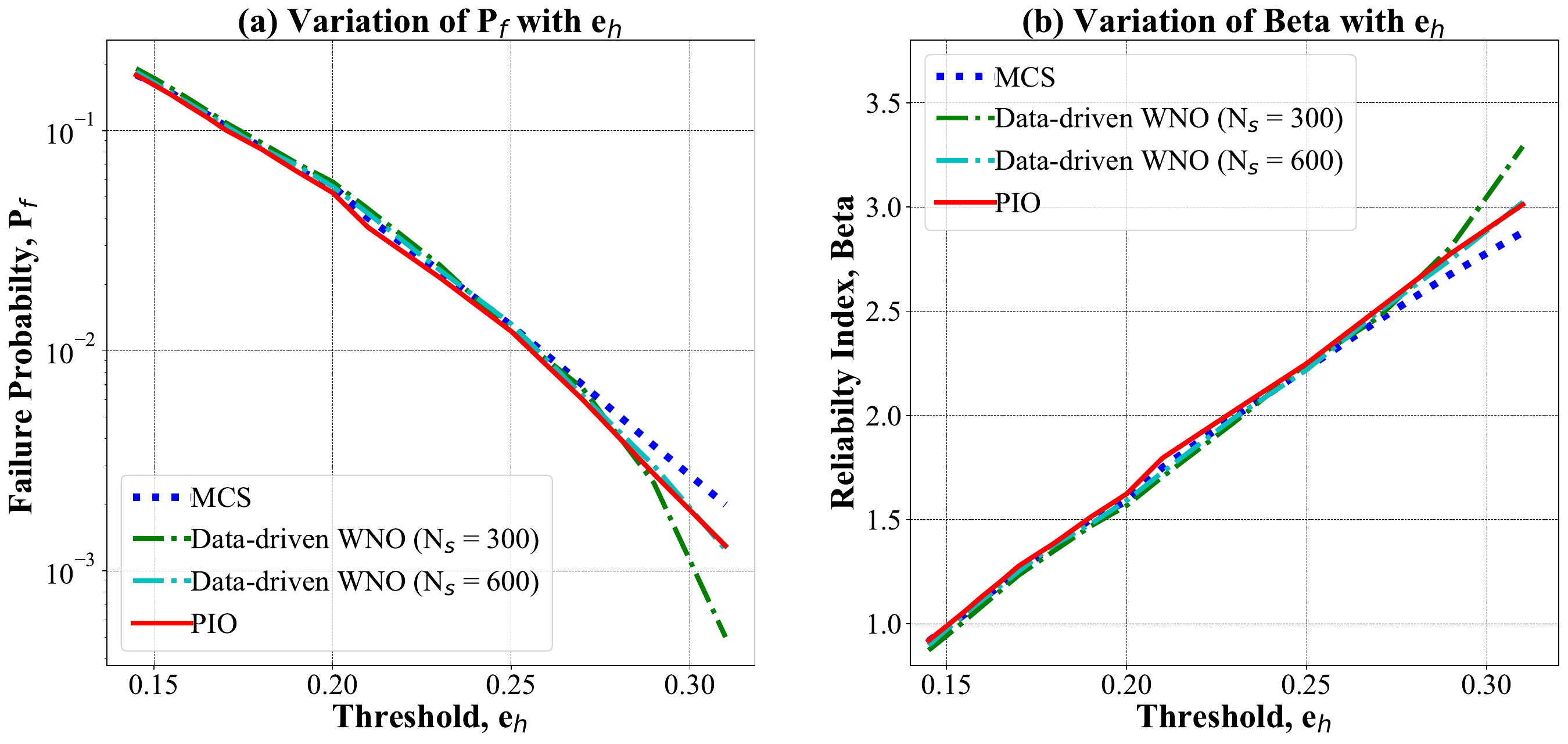}}
    \caption{Prediction results of the example: impulse transmission in nerve for the failure probability ($P_f$) and the reliability index ($\beta$) with unceasing limit state threshold, obtained by PIO in comparison with the results of MCS, and data-driven WNOs (trained with a number of samples, $N_s = 300$ and $N_s = 600$)}
    \label{fig:nagumo_vare_h}
\end{figure}

Firstly, we validate the efficacy of the PIO-predicted response against the ground truths obtained from the numerical solvers. \autoref{fig:nagumophy} shows a visualisation of the predictions in comparison with ground truths, demonstrating excellent agreement. For reliability analysis, the grid and the threshold in the limit-state function are specified to be $x_{sp} = 32$ and $e_h = 1.45$. Based on the limit state, a probability density plot of the first passage failure time is shown in \autoref{fig:nagumopdf}. Case studies by varying the number of training inputs are illustrated. We observe that the result continuously improves with an increase in the input samples. The fact that the proposed approach accurately captures the multimodal distribution is impressive and indicative of the strength of the model.
It is again important to note that PIO is trained only using initial conditions, and hence, there is no need of generating response data by running computationally expensive simulations. 
The failure probability and the reliability index ($\beta$) with varying numbers of training input samples are also provided in the \autoref{table2}. It is observed that the failure probability and reliability index obtained using the proposed approach match well with the benchmark solution generated using MCS. 
Variation of failure probability and reliability indices with change in threshold is shown in \autoref{fig:nagumo_vare_h}. 
It is observed that both PIO and data-driven WNO with 600 training samples yield excellent results that closely follows the ground truth obtained using MCS. It is important to note that, unlike PIO, data-driven WNO necessitates running a computationally expensive simulator to generate training outputs and, hence, is significantly more expensive from a computational point-of-view. Given that the training and inference time for PIO and data-driven WNO are similar, it can be concluded that the data-driven WNO requires $600\times t_s$ additional time where $t_s$ is the time required for running a single simulation. Data-driven WNO with 300 training samples requires yields erroneous results for thresholds corresponding to small failure probabilities.

\begin{table}[ht!]
    \centering
    \caption{Results of first passage failure probability obtained using the proposed framework for the case of the impulse transmission in nerve example with varying trained initial conditions.}
    \label{table2}
\begin{tabular}{lccccc} 
\hline
\textbf{No. of training conditions}& \textbf{300} & \textbf{400} & \textbf{500}&\textbf{600}&\textbf{Actual}\\
\hline
Failure probability ($P_f$) & 0.150 & 0.170 & $0.174$ & 0.177 & 0.1795\\ \hline
Reliability index ($\beta$) & 1.034 &  0.954  &  0.940 & 0.927	& 0.917	\\ 
\hline 
\end{tabular}
\end{table}

\subsection{Fluid flow through a porous medium}\label{exampl3}
\begin{figure}[!h]
    \centering{
    \includegraphics[width=0.4 \textwidth]{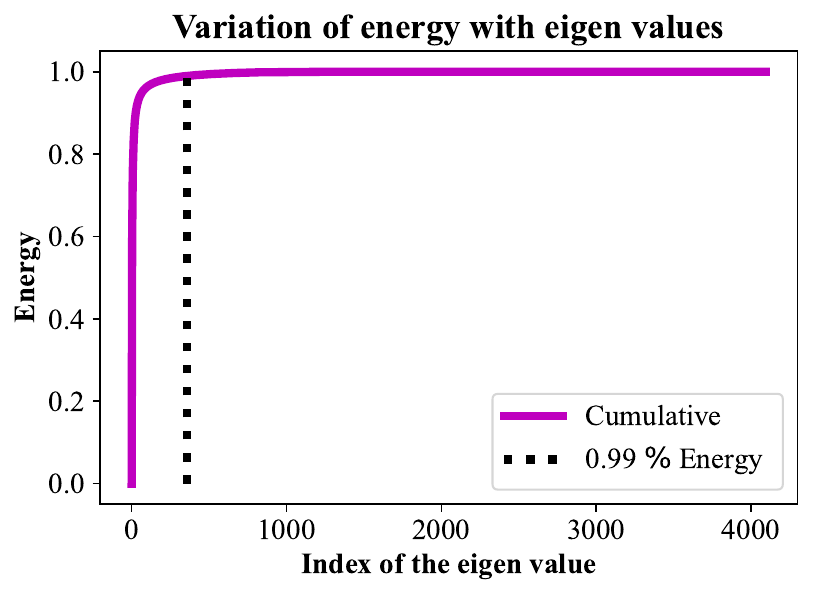}}
    \caption{Intrinsic dimensionality of the input function obtained for the example of flow through a porous medium, where the 99$\%$ of the energy is contained in the first 358 eigenvalues.}
    \label{fig:Darcys_Intrinsic}
\end{figure}

Flow through a porous medium is a well-studied phenomenon in the petroleum industry. Understanding the flow of fluids through reservoir rocks is crucial, and this process can be effectively described using Darcy's flow equation \cite{dejam2017pre}. The model can be used to estimate pressure difference, which allows a certain volume of fluid with a specific viscosity to flow through a given cross-sectional area in a unit of time under a given permeability of a rock. For a stationary flow, the governing PDE for Darcy flow is expressed as:
\begin{equation}
\begin{aligned}
-\nabla \cdot (a(x, y) \nabla u(x, y)) & = f(x, y), \quad x, y \in (0, R)
u(x, y) & = u_0(x, y), \quad x, y \in \partial(0, R)
\end{aligned}
\end{equation}
where $u(x, y) = u_0(x, y)$ is the boundary condition satisfying the zero Dirichlet boundary condition ($u_0(x, y) = 0$),  $ a(x, y)$ is the permeability field and $u(x, y)$ is the corresponding pressure field. $f(x, y)$ represents a source function, and it is taken to be a constant value of 1 ($f(x,y)=1$). The flow problem is confined to the domain of unit square size with $x, y \in (0, 1)^2$. To analyse the reliability of the flow, it is important to monitor variations in pressure, where an increase in pressure value beyond a certain limit can cause uncontrollable flow, leading to potential blowouts or structural damage to the reservoir. Taking this into consideration, the limit state function is defined here as a maximum absolute value of the pressure field ($u$) that crosses threshold $e_h$. It is mathematically expressed as:
\begin{equation}
    \mathcal{J}(x)=e_h-{|u(x,y)|}_{max},
\end{equation}
where the threshold value for the pressure output field is set to be $e_h = 0.078$. To carry out the reliability analysis, we consider the permeability field to be a random field generated from a distribution such that $a = N(0,(-\Delta + 9I)^{-2})$ with zero Neumann boundary conditions on the Laplacian. The input resolution is chosen to be $64 \times 64$, and the generated random field has an intrinsic dimensionality of 358 (\autoref{fig:Darcys_Intrinsic}). This indicates the high dimensionality of the input field. The PIO is trained to learn the mapping between the stochastic permeability field and the pressure field, $\mathcal{M}: a(x, y) \mapsto u(x, y)$. Similar to previous examples, the operator is trained with different sizes to training inputs to illustrate its convergence.

\begin{figure}[!h]
    \centering{
    \includegraphics[width=1\textwidth]{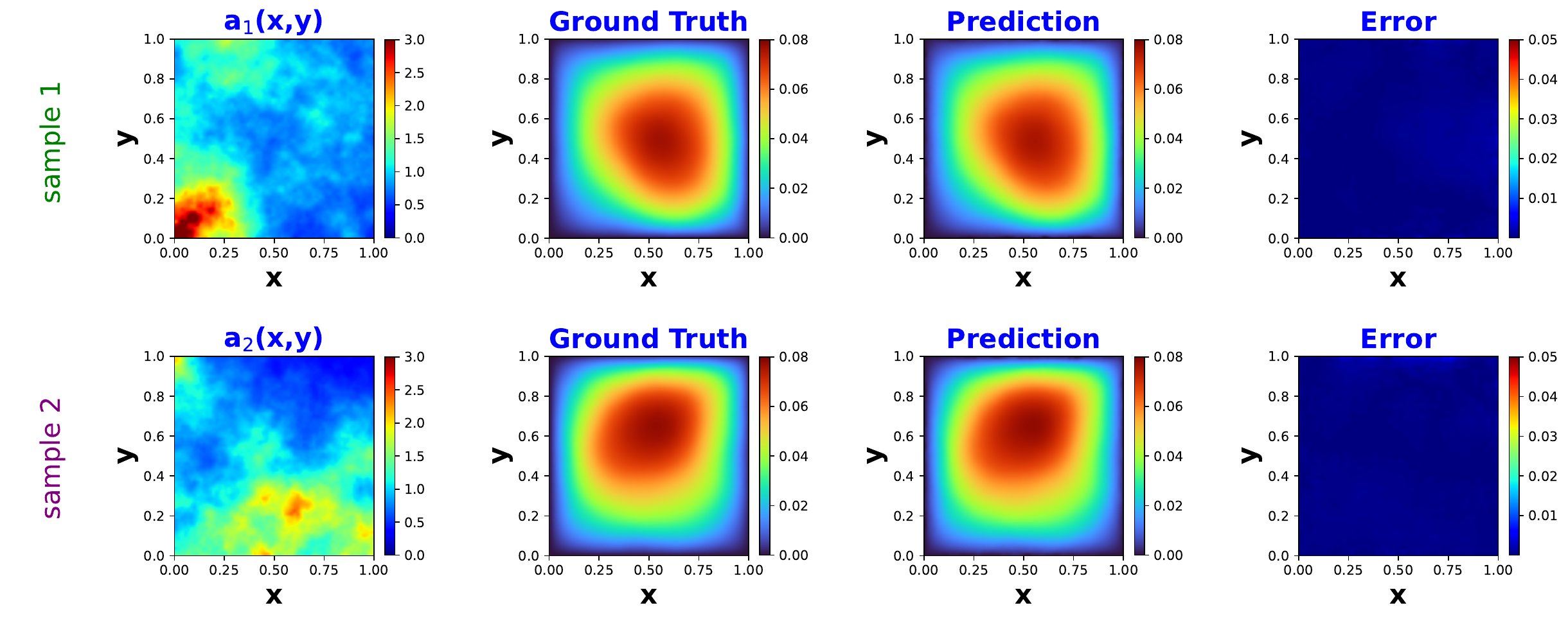}}
    \caption{The results for the example of flow through a porous medium comprised parametric permeability fields, ground truth solutions, predictions, and error plots demonstrated for using  2 different unseen sample instances. The PIO effectively maps the permeability fields to the corresponding pressure field $u(x,y)$ over the domain, with a spatial resolution of $64 \times 64$.}
    \label{fig:darcyphy}
\end{figure}

\begin{figure}[!h]
    \centering{
    \includegraphics[width=0.8\textwidth]{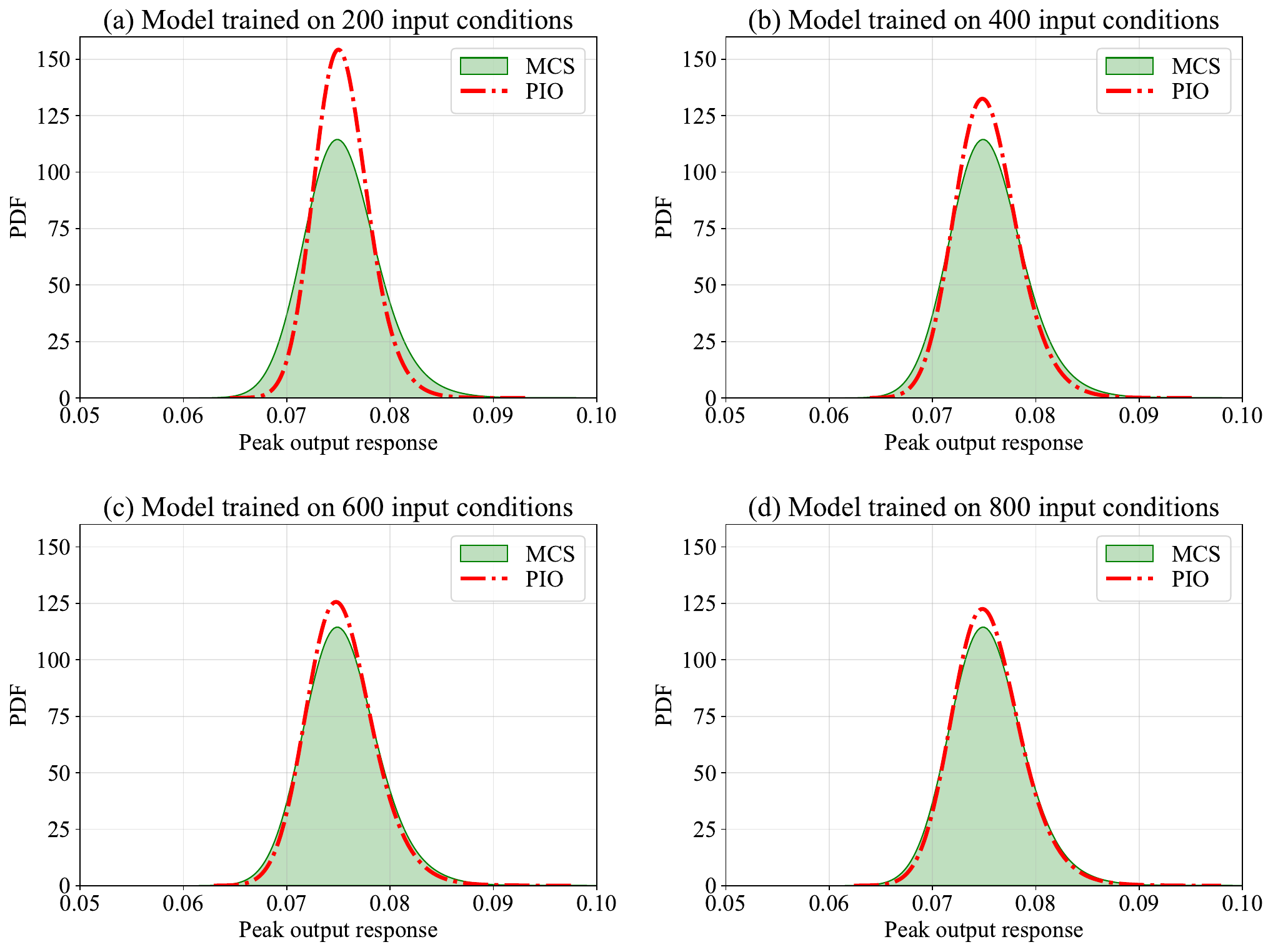}}
    \caption{PDF plots of peak pressure values obtained by MCS and PIO trained with varying number of permeability field conditions for the example of flow through a porous medium}
    \label{fig:darcypdf}
\end{figure}

\begin{figure}[!h]
    \centering{
    \includegraphics[width=0.8\textwidth]{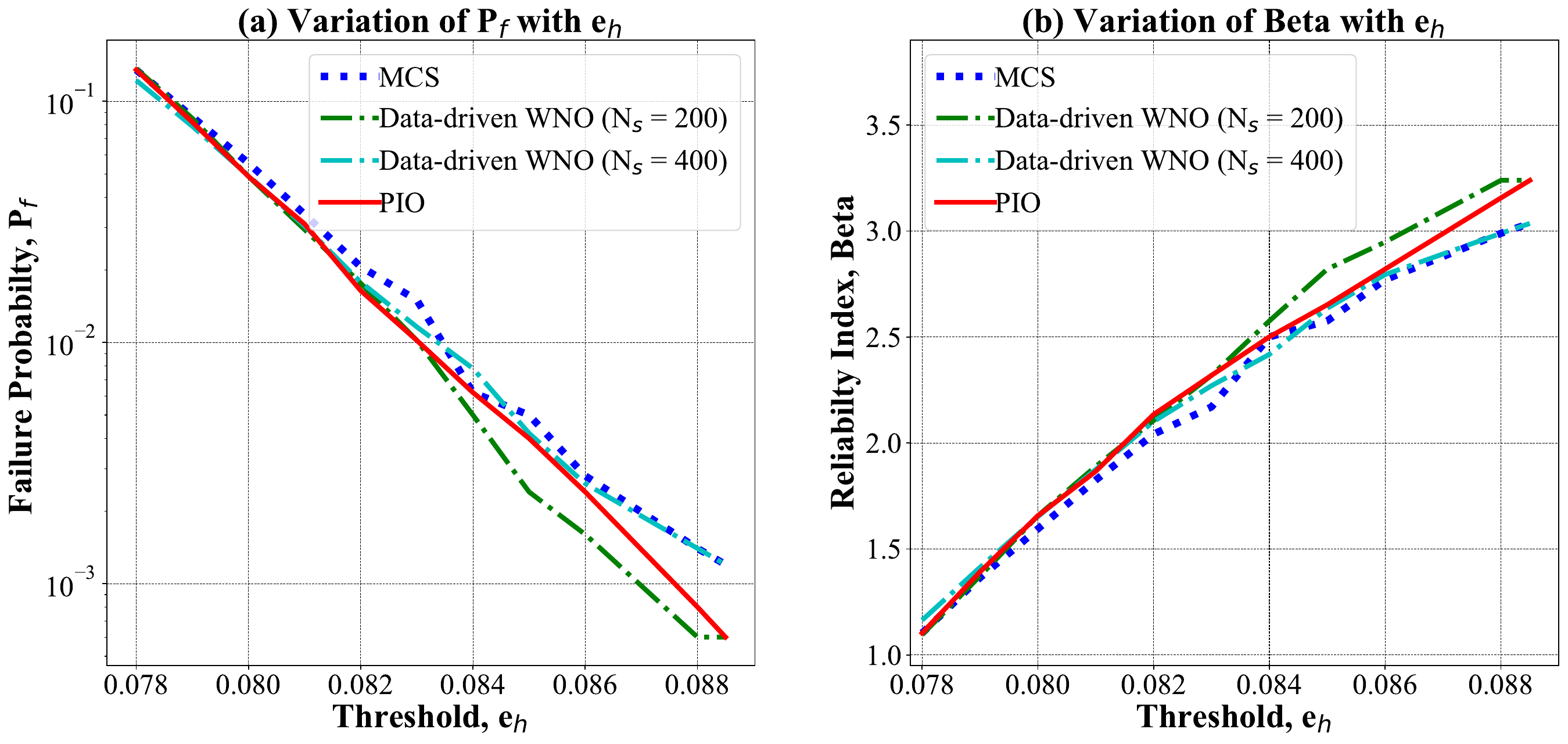}}
    \caption{Prediction results of the example: fluid flow through a porous medium for the failure probability ($P_f$) and the reliability index ($\beta$) with unceasing limit state threshold, obtained by PIO in comparison with the results of MCS, and data-driven WNOs (trained with a number of samples, $N_s = 200$ and $N_s = 400$)}
    \label{fig:darcy_var_eh}
\end{figure}

Before proceeding with the reliability analysis, the predictions obtained using PIO are validated with the numerical solutions. The results are presented in \autoref{fig:darcyphy} indicate that PIO yields highly accurate results. Moving on to detailed results of reliability estimation, it is noted here that contrary to the first two examples, the third example presented here is a time-independent problem. Thus, to evaluate the efficacy, instead of the PDF plot of the failure time, a probability density plot of the maximum pressure value across the samples is illustrated in \autoref{fig:darcypdf}. As expected, we observe a gradual convergence with increase in the number of training inputs. 
The failure probability reliability index obtained using different training sample sizes are shown in \autoref{table3}. The probability of failure and reliability index obtained with 800 training inputs match exactly with the benchmark results obtained using MCS. Lastly, a comparative study on the reliability estimates corresponding to different thresholds is illustrated in \autoref{fig:darcy_var_eh}. We observe that PIO yields accurate results until $e_h = 0.086$. Data-driven WNO with 200 training samples yields accurate results until $e_h = 0.083$. Data-driven WNO with 400 training samples yields the best results and yields accurate results throughout the threshold limits illustrated in the figure. However, it is to be noted that data-driven WNO requires simulation data and hence, the overall computational cost is directly proportional to the number of training samples. Therefore, for this example, PIO can be viewed as an efficient alternative to the data-driven WNO with 400 training samples. Note that PIO also provides the option of using both data and physics; however, the same is not explored here as the objective is to investigate the capability of PIO to solve reliability analysis problems from no simulation data.



\begin{table}[ht!]
    \centering
    \caption{Results of failure probability obtained using   
    proposed framework for the example of flow through a porous medium with varying trained permeability fields.}
    \label{table3}
\begin{tabular}{lccccc} 
\hline
\textbf{No. of training conditions}& \textbf{200} & \textbf{400} & \textbf{600}&\textbf{800}&\textbf{Actual}\\
\hline
Failure probability ($P_f$) &  0.0948  &  0.1146 &   $0.1214$ & 0.1350 & 0.1350\\ \hline
Reliability index ($\beta$) &  1.3118 &  1.202  &  1.168 & 1.103 & 1.103\\ 
\hline 
\end{tabular}
\end{table}
\subsection{Phase transitions in alloys}\label{example4}
\begin{figure}[!h]
    \centering{
    \includegraphics[width=0.4 \textwidth]{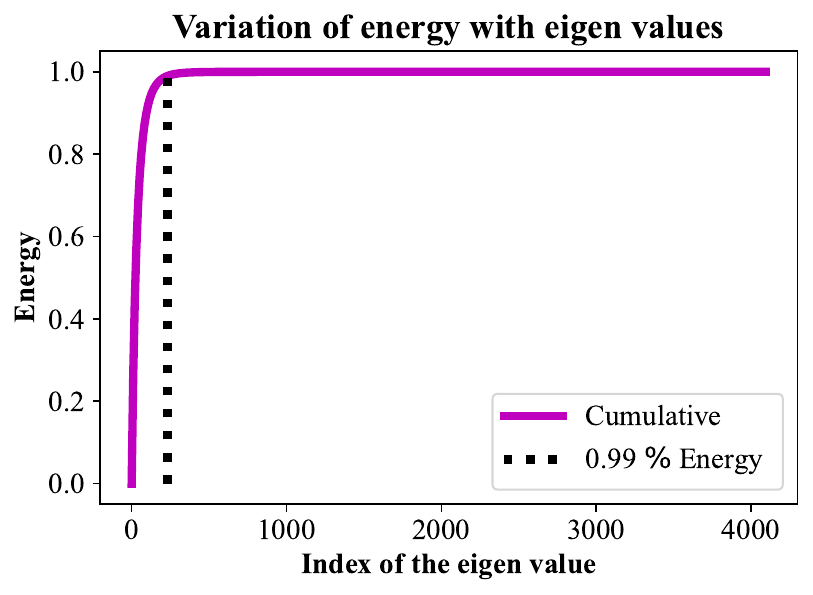}}
    \caption{Intrinsic dimensionality of the input function obtained for example of phase transitions in alloys, where the 99$\%$ of the energy is contained in the first 232 eigenvalues.}
    \label{fig:Allencahn_Intrinsic}
\end{figure}

In the last example, we examine the reliability of the phase transition phenomenon in alloys. Analysing the phase transitions in alloys is crucial for understanding and controlling the properties of materials used in various industrial applications. Modeling these transitions provides insights into the microstructural evolution of alloys, which directly affects their mechanical, thermal, and electrical properties. Effective modelling enables us to predict the behaviour of alloys under different conditions and optimise material properties for specific applications. Thus, proactive steps can be taken to avoid design failures stemming from the material properties. The dynamics of phase separation in a microstructure formation can be effectively modelled using the Allen-Cahn equation \cite{hussain2019approximate,ma2017numerical}. The equation describes the evolution of the order parameter $u$, which represents the different phases in the alloy over time and space.  An increase in the order parameter above certain limits can lead to significant microstructural changes, potentially altering the material properties and impacting the performance of the alloy. Thus, the limit state in this context can be defined by the maximum value over the 2-D spatiotemporal order parameter field $u$. Specifically, the failure is considered if the maximum absolute value of the output response $u$ between the time steps, $[t_1,t_2]$ exceeds the threshold $e_h$. Mathematically, this can be expressed as:
\begin{equation}\label{limit4}
    \mathcal{J}(x) = e_h - {|u(x,y,t)|}_{max},
\end{equation}
where $u(x,y,t)$ represents the solution of the two-dimensional time-dependent Allen–Cahn equation is formulated as:
\begin{equation}
\begin{aligned}
\partial_t u(x, y, t) & =\epsilon \Delta u(x, y, t)+u(x, y, t)-u(x, y, t)^3, & x, y \in(0,1) \\
u(x, y, 0) & =u_0(x, y) & x, y \in(0,1)
\end{aligned}
\end{equation}
In the above expression, $\epsilon$ denotes the viscosity coefficient such that $\epsilon \in \mathbb{R}^{+*}$ (positive real) , which is set to $\epsilon = 1 \times 10^{-3}$ in this study. The problem is defined on a periodic boundary. For the reliability analysis of the problem, we consider the initial conditions to be random fields and generated using a Gaussian random field with the following kernel:
\begin{equation}
    \mathcal{K}(x, y)=\tau^{(\alpha-1)}\left(\pi^2\left(x^2+y^2\right)+\tau^2\right)^{\frac{\pi}{2}},
\end{equation}
where, the parameters for the kernel are fixed to \(\tau = 15\) and \(\alpha = 1\). The intrinsic dimension of the input field is calculated as 232 based on the estimation given in \autoref{fig:Allencahn_Intrinsic}. Therefore, this can be categorized as a high-dimensional time-dependent reliability analysis problem.
We train the PIO to learn the mapping 

Our objective is to first obtain the predictions of temporal responses using the solution operator. The operator in this problem is trained to map $\bm{u}$ from the domain $(0,1)^2 \times [0,10]$ to $(0,1)^2 \times (10, T]$. For training the PIO, 600 different input conditions are used, where the resolution is chosen to be $64 \times 64$. For illustration, we set the target time steps to be $T=22$, i.e., we predict the solution for the next 13 steps. 

\begin{figure}[!h]
    \centering{
    \includegraphics[width=1\textwidth]{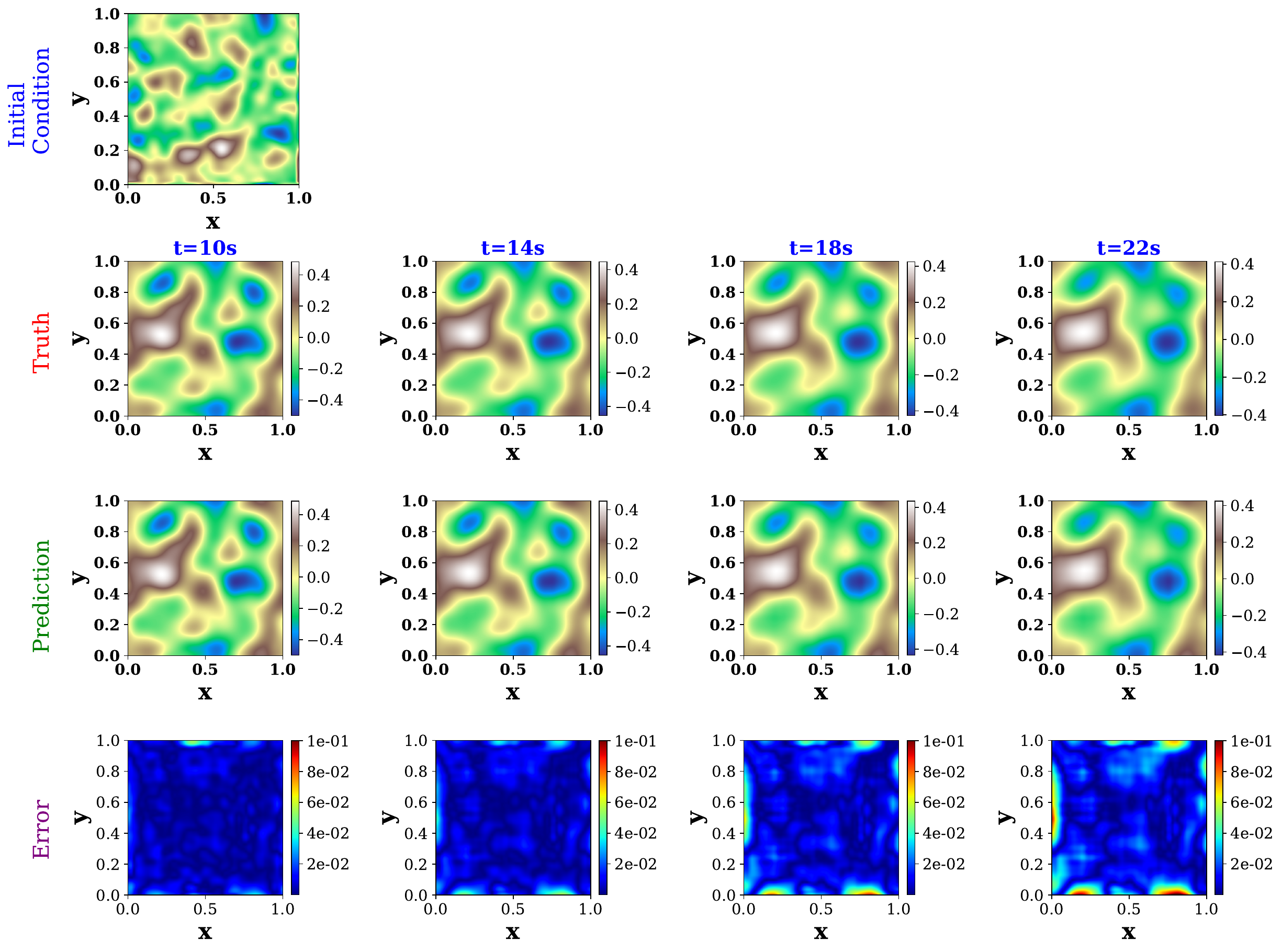}}
    \caption{The results for example of phase transitions in alloys comprised of given initial field and corresponding ground truth solutions, predictions, and error plots illustrated with an unseen sample instance at time steps $10~s,14~s,18~s$ and $22~s$. The PIO receives the spatial field, $u(x,y)$ with a resolution of $65\times65$, for the initial 10-time steps and which maps to the corresponding solution $u(x,y)$ for the 13 future time-steps}
    \label{fig:Allencahn_physics}
\end{figure}

\begin{figure}[!h]
    \centering{
    \includegraphics[width=1\textwidth]{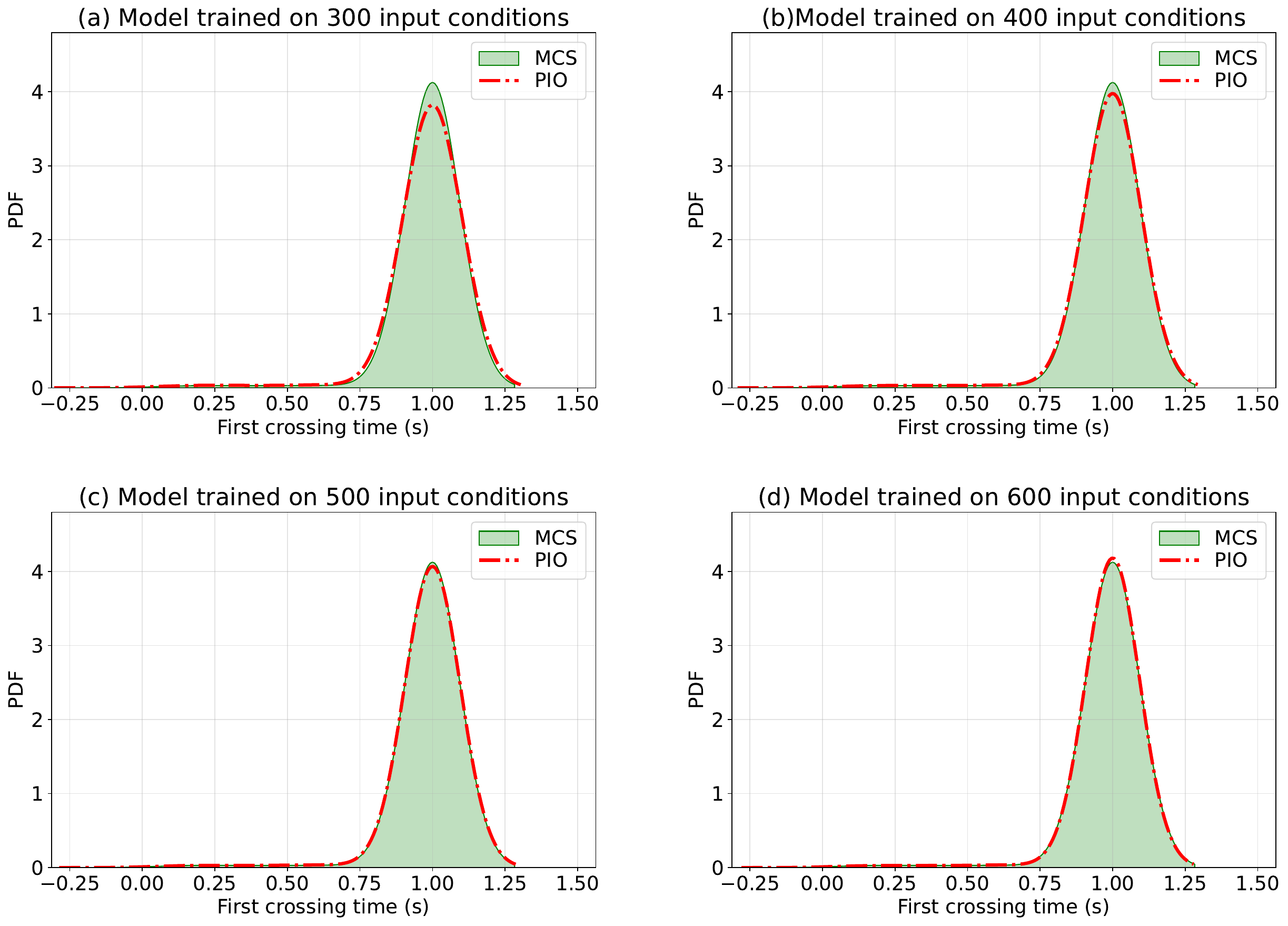}}
    \caption{PDF plots of failure time obtained by MCS and PIO for varying limiting state conditions for example of phase transitions in alloys}
    \label{fig:Allencahn_pdf}
\end{figure}

\begin{figure}[!h]
    \centering{
    \includegraphics[width=0.8\textwidth]{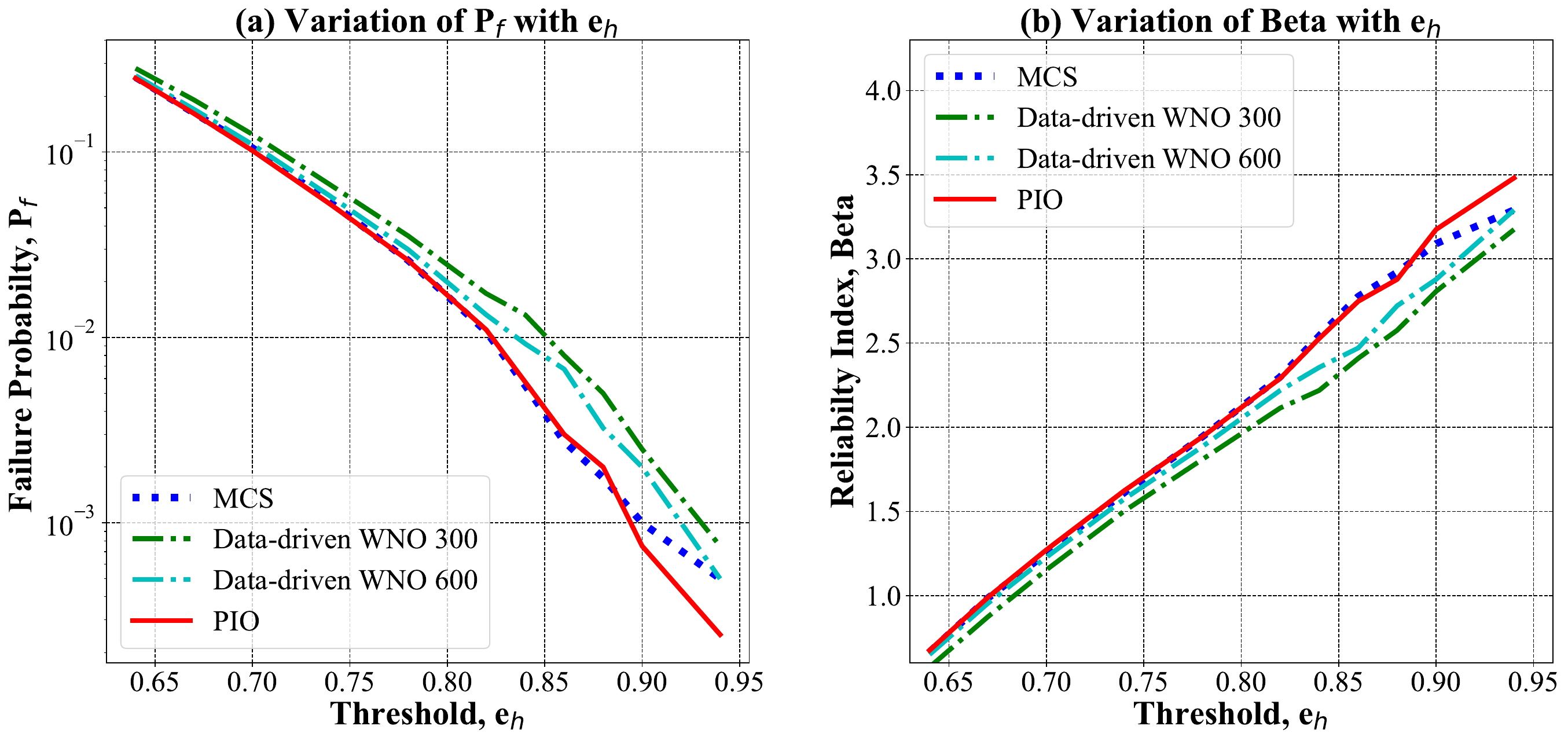}}
    \caption{Prediction results of the example: phase transitions in alloys for the failure probability ($P_f$) and the reliability index ($\beta$) with unceasing limit state threshold, obtained by PIO in comparison with the results of MCS, and data-driven WNOs (trained with a number of samples, $N_s = 300$ and $N_s = 600$)}
    \label{fig:allencahn_var_eh}
\end{figure}

Firstly, the prediction results of the PIO are compared with the ground truth solutions to validate the efficacy of the operator. The results are presented in \autoref{fig:Allencahn_physics}. It is evident from the results that the predictions of the PIO almost exactly emulate the ground truth solutions. To obtain the PDF plot and to compute reliability, we employ the above-mentioned limit state function given in \autoref{limit4}. For estimating failure probability, output responses between the time step $10-22$, $([t_1, t_2])$ are considered.  The value of the threshold is chosen such that $e_h=0.78$. 
The results of the study, where PIO is trained with different numbers of input instances, are shown in the \autoref{fig:Allencahn_pdf}. Subsequently, \autoref{table4} lists the estimated failure probabilities and reliability indices. The results presented in the table and the PDF plots demonstrate the efficacy of the PIO in estimating failure probability. It also reinforces the fact that the accuracy of the reliability estimation is increased with an increase in the number of initial training inputs. A comparative study, which includes the results obtained by the PIO, data-driven WNO and the MCS for a set of varied limit state thresholds, is presented in \autoref{fig:allencahn_var_eh}. We observe that PIO yields the best results across different thresholds. Data-driven WNO, on the other hand, fails to capture the failure probability beyond $e_h=0.80$. This illustrates the benefit of PIO for systems where the underlying evolution law is relatively complex. Overall, the fact that PIO can solve time-dependent reliability analysis problems of such high-dimensional systems from no simulation data is impressive and indicates its potential for the actual application of this method.

\begin{table}[ht!]
    \centering
    \caption{Results of first passage failure probability obtained using   
    proposed framework for the example of phase transitions in alloys with varying trained initial conditions.}
    \label{table4}
\begin{tabular}{lccccc} 
\hline
\textbf{No. of training conditions}& \textbf{300} & \textbf{400} & \textbf{500}&\textbf{600}&\textbf{Actual}\\
\hline
Failure probability ($P_f$) & 0.03075 & 0.02825 & $0.02725$ & 0.026 & 0.026\\ \hline
Reliability index ($\beta$) & 1.8699 &  1.9072  &  1.9228 & 1.943	& 1.943	\\ 
\hline 
\end{tabular}
\end{table}

\section{Conclusions}\label{sec:Conclusions}
Operator learning is a relatively new paradigm in deep learning that aims to learn the mapping between function spaces. Recent studies have illustrated the huge potential offered by operator learning algorithms.
In this work, we investigated the possible application of Physics-Informed Operator (PIO) learning paradigm for solving reliability analysis problems. In PIO, physics-informed loss function formulated using the governing physics is incorporated within the operator learning algorithm to alleviate the data-hungry nature of data-driven operator learning algorithms. This is of particular interest in reliability engineering, as one of the major goals in reliability analysis is to minimize the computational cost by minimizing the number of simulations. However, solving reliability analysis problems is significantly more challenging than quantifying the response statistics as the tails region of the probability density function is to be captured to obtain accurate estimates of the reliability. Accordingly, in this work, we investigate the possibility of using PIO for solving reliability analysis problems. We investigate the applicability of PIO for both time-dependent and time-independent reliability analysis problems for four canonical systems governed by partial differential equations. The major findings of this investigation are summarized below:
\begin{itemize}
    \item \textbf{Data efficiency}: We observed that PIO trained only using the physics loss, boundary conditions, and initial conditions yields satisfactory results for all the examples. This indicates that PIO can eliminate the requirement of performing simulation to quantify the uncertainty. In other words, PIO  is significantly for efficient than existing data-driven surrogate models used for reliability analysis.
    \item \textbf{Scalability}: The application of physics-informed loss to eliminate computational efficiency is not new, as the same is also offered by physics informed neural networks. However, physics-informed neural networks are not scalable. PIO, on the other hand, is scalable, as illustrated in this work. For instance, the flow through porous media and the phase transition in alloy problems solved in this paper involved 358 and 243 intrinsic dimensionality (random variables), respectively.
    \item \textbf{Accuracy}: For all the examples presented, PIO yields reasonably accurate estimates of the probability of failure at different thresholds. In general, for time-dependent reliability analysis problems, PIO either yield comparable (for the first two examples) or better results (for the last example) when compared with the data-driven counterparts. A key observation was that PIO outperforms the data-driven counterpart as the complexity of the underlying evolution equation (governing physics) increases. This was evident from the phase transition in metals example presented. However, for time-independent reliability analysis, data-driven operator learning outperformed the PIO for relatively small failure probabilities.
\end{itemize}
Based on these observations, the potential offered by PIO is evident. Having said that, it may be noted that this study only investigates the potential application of physics-informed wavelet neural operator. However, there exist other physics-informed operators, including physics-informed DeepONet, that can also be used. Given the comparable performance of the data-driven variants of these operators, we expect similar observations to hold for the physics-informed counterparts of these operators. Also, throughout this study, we have assumed that the exact physics is known. In case the exact physics is not known, direct application of PIO will not be possible. We have addressed this in a separate study.

\section*{Acknowledgements}
NN acknowledge the support received from the Ministry of Education in the form of the Prime Minister's Research Fellowship, and T acknowledges the support received from the Ministry of Education through a Senior Research Fellowship. SC acknowledges the financial support received from the Ministry of Ports and Shipping via letter number ST-14011/74/MT (356529) and Science and Engineering Research Board via grant number CRG/2023/007667

% \bibliographystyle{ieeetr}
% \bibliography{references}  %%% Uncomment this line and comment out the ``thebibliography'' section below to use the external .bib file (using bibtex).

\begin{thebibliography}{10}

\bibitem{ni2020reliability}
P.~Ni, J.~Li, H.~Hao, W.~Yan, X.~Du, H.~Zhou, Reliability analysis and design optimization of nonlinear structures, Reliability engineering \& system safety 198 (2020) 106860.

\bibitem{matteo2021time}
F.~Matteo, G.~Carlo, P.~Federico, Z.~Enrico, et~al., Time-dependent reliability analysis of the reactor building of a nuclear power plant for accounting of its aging and degradation, Reliability Engineering \& System Safety 205 (2021) 107173.

\bibitem{de2005naroas}
A.~G. de~Araujo~G{\'o}es, M.~A.~B. Alvarenga, P.~F. e~Melo, Naroas: a neural network-based advanced operator support system for the assessment of systems reliability, Reliability Engineering \& System Safety 87~(2) (2005) 149--161.

\bibitem{tao2024reliability}
H.~Tao, P.~Jia, X.~Wang, L.~Wang, Reliability analysis of subsea control module based on dynamic bayesian network and digital twin, Reliability Engineering \& System Safety 248 (2024) 110153.

\bibitem{thakur1978monte}
R.~Thakur, K.~Misra, Monte carlo simulation for reliability evaluation of complex systems, International Journal of Systems Science 9~(11) (1978) 1303--1308.

\bibitem{rubinstein2016simulation}
R.~Y. Rubinstein, D.~P. Kroese, Simulation and the Monte Carlo method, Vol.~10, John Wiley \& Sons, 2016.

\bibitem{xiong2021fast}
Y.~Xiong, S.~Sampath, A fast-convergence algorithm for reliability analysis based on the ak-mcs, Reliability Engineering \& System Safety 213 (2021) 107693.

\bibitem{chang2022mc}
P.-C. Chang, Mc-based simulation approach for two-terminal multi-state network reliability evaluation without knowing d-mcs, Reliability Engineering \& System Safety 220 (2022) 108289.

\bibitem{au1999new}
S.-K. Au, J.~L. Beck, A new adaptive importance sampling scheme for reliability calculations, Structural safety 21~(2) (1999) 135--158.

\bibitem{li2005curse}
B.~Li, T.~Bengtsson, P.~Bickel, Curse-of-dimensionality revisited: Collapse of importance sampling in very large scale systems, Rapport technique 85 (2005) 205.

\bibitem{engelund1993benchmark}
S.~Engelund, R.~Rackwitz, A benchmark study on importance sampling techniques in structural reliability, Structural safety 12~(4) (1993) 255--276.

\bibitem{au2001estimation}
S.-K. Au, J.~L. Beck, Estimation of small failure probabilities in high dimensions by subset simulation, Probabilistic engineering mechanics 16~(4) (2001) 263--277.

\bibitem{au2014engineering}
S.-K. Au, Y.~Wang, Engineering risk assessment with subset simulation, John Wiley \& Sons, 2014.

\bibitem{zuev2015subset}
K.~Zuev, Subset simulation method for rare event estimation: an introduction, arXiv preprint arXiv:1505.03506 (2015).

\bibitem{ditlevsen1990general}
O.~Ditlevsen, R.~E. Melchers, H.~Gluver, General multi-dimensional probability integration by directional simulation, Computers \& Structures 36~(2) (1990) 355--368.

\bibitem{xu2021machine}
Z.~Xu, J.~H. Saleh, Machine learning for reliability engineering and safety applications: Review of current status and future opportunities, Reliability Engineering \& System Safety 211 (2021) 107530.

\bibitem{blatman2011adaptive}
G.~Blatman, B.~Sudret, Adaptive sparse polynomial chaos expansion based on least angle regression, Journal of computational Physics 230~(6) (2011) 2345--2367.

\bibitem{sudret2008global}
B.~Sudret, Global sensitivity analysis using polynomial chaos expansions, Reliability engineering \& system safety 93~(7) (2008) 964--979.

\bibitem{de2013new}
S.~De~Marchi, G.~Santin, A new stable basis for radial basis function interpolation, Journal of Computational and Applied Mathematics 253 (2013) 1--13.

\bibitem{li2018sequential}
X.~Li, C.~Gong, L.~Gu, W.~Gao, Z.~Jing, H.~Su, A sequential surrogate method for reliability analysis based on radial basis function, Structural Safety 73 (2018) 42--53.

\bibitem{ziehn2009gui}
T.~Ziehn, A.~S. Tomlin, Gui--hdmr--a software tool for global sensitivity analysis of complex models, Environmental Modelling \& Software 24~(7) (2009) 775--785.

\bibitem{metya2017system}
S.~Metya, T.~Mukhopadhyay, S.~Adhikari, G.~Bhattacharya, System reliability analysis of soil slopes with general slip surfaces using multivariate adaptive regression splines, Computers and Geotechnics 87 (2017) 212--228.

\bibitem{bilionis2012multi}
I.~Bilionis, N.~Zabaras, Multi-output local gaussian process regression: Applications to uncertainty quantification, Journal of Computational Physics 231~(17) (2012) 5718--5746.

\bibitem{bilionis2013multi}
I.~Bilionis, N.~Zabaras, B.~A. Konomi, G.~Lin, Multi-output separable gaussian process: Towards an efficient, fully bayesian paradigm for uncertainty quantification, Journal of Computational Physics 241 (2013) 212--239.

\bibitem{tripathy2016gaussian}
R.~Tripathy, I.~Bilionis, M.~Gonzalez, Gaussian processes with built-in dimensionality reduction: Applications to high-dimensional uncertainty propagation, Journal of Computational Physics 321 (2016) 191--223.

\bibitem{atkinson2018structured}
S.~Atkinson, N.~Zabaras, Structured bayesian gaussian process latent variable model, arXiv preprint arXiv:1805.08665 (2018).

\bibitem{atkinson2019structured}
S.~Atkinson, N.~Zabaras, Structured bayesian gaussian process latent variable model: Applications to data-driven dimensionality reduction and high-dimensional inversion, Journal of Computational Physics 383 (2019) 166--195.

\bibitem{bansal2022physics}
P.~Bansal, Z.~Zheng, C.~Shao, J.~Li, M.~Banu, B.~E. Carlson, Y.~Li, Physics-informed machine learning assisted uncertainty quantification for the corrosion of dissimilar material joints, Reliability Engineering \& System Safety 227 (2022) 108711.

\bibitem{qian2021time}
H.-M. Qian, Y.-F. Li, H.-Z. Huang, Time-variant system reliability analysis method for a small failure probability problem, Reliability Engineering \& System Safety 205 (2021) 107261.

\bibitem{ling2020efficient}
C.~Ling, Z.~Lu, X.~Zhang, An efficient method based on ak-mcs for estimating failure probability function, Reliability Engineering \& System Safety 201 (2020) 106975.

\bibitem{roy2020support}
A.~Roy, S.~Chakraborty, Support vector regression based metamodel by sequential adaptive sampling for reliability analysis of structures, Reliability Engineering \& System Safety 200 (2020) 106948.

\bibitem{ghosh2018support}
S.~Ghosh, A.~Roy, S.~Chakraborty, Support vector regression based metamodeling for seismic reliability analysis of structures, Applied Mathematical Modelling 64 (2018) 584--602.

\bibitem{cheng2021adaptive}
K.~Cheng, Z.~Lu, Adaptive bayesian support vector regression model for structural reliability analysis, Reliability Engineering \& System Safety 206 (2021) 107286.

\bibitem{zhou2022towards}
T.~Zhou, T.~Han, E.~L. Droguett, Towards trustworthy machine fault diagnosis: A probabilistic bayesian deep learning framework, Reliability Engineering \& System Safety 224 (2022) 108525.

\bibitem{pepper2022adaptive}
N.~Pepper, L.~Crespo, F.~Montomoli, Adaptive learning for reliability analysis using support vector machines, Reliability Engineering \& System Safety 226 (2022) 108635.

\bibitem{chen2022support}
J.-Y. Chen, Y.-W. Feng, D.~Teng, C.~Lu, C.-W. Fei, Support vector machine-based similarity selection method for structural transient reliability analysis, Reliability Engineering \& System Safety 223 (2022) 108513.

\bibitem{elhewy2006reliability}
A.~H. Elhewy, E.~Mesbahi, Y.~Pu, Reliability analysis of structures using neural network method, Probabilistic Engineering Mechanics 21~(1) (2006) 44--53.

\bibitem{navaneeth2022koopman}
N.~Navaneeth, S.~Chakraborty, Koopman operator for time-dependent reliability analysis, Probabilistic Engineering Mechanics 70 (2022) 103372.

\bibitem{navaneeth2022surrogate}
N.~Navaneeth, S.~Chakraborty, Surrogate assisted active subspace and active subspace assisted surrogate—a new paradigm for high dimensional structural reliability analysis, Computer Methods in Applied Mechanics and Engineering 389 (2022) 114374.

\bibitem{chakraborty2017efficient}
S.~Chakraborty, R.~Chowdhury, An efficient algorithm for building locally refined hp--adaptive h-pcfe: Application to uncertainty quantification, Journal of Computational Physics 351 (2017) 59--79.

\bibitem{chakraborty2017hybrid}
S.~Chakraborty, R.~Chowdhury, Hybrid framework for the estimation of rare failure event probability, Journal of Engineering Mechanics 143~(5) (2017) 04017010.

\bibitem{wu2023adaptive}
H.~Wu, Y.~Xu, Z.~Liu, Y.~Li, P.~Wang, Adaptive machine learning with physics-based simulations for mean time to failure prediction of engineering systems, Reliability Engineering \& System Safety 240 (2023) 109553.

\bibitem{mathpati2023mantra}
Y.~C. Mathpati, K.~S. More, T.~Tripura, R.~Nayek, S.~Chakraborty, Mantra: A framework for model agnostic reliability analysis, Reliability Engineering \& System Safety 235 (2023) 109233.

\bibitem{chakraborty2023deep}
S.~Chakraborty, et~al., Deep physics corrector: A physics enhanced deep learning architecture for solving stochastic differential equations, Journal of Computational Physics 479 (2023) 112004.

\bibitem{chakraborty2023dpa}
S.~Chakraborty, et~al., Dpa-wno: A gray box model for a class of stochastic mechanics problem, arXiv preprint arXiv:2309.15128 (2023).

\bibitem{ramabathiran2021spinn}
A.~A. Ramabathiran, P.~Ramachandran, Spinn: Sparse, physics-based, and partially interpretable neural networks for pdes, Journal of Computational Physics 445 (2021) 110600.

\bibitem{sirignano2018dgm}
J.~Sirignano, K.~Spiliopoulos, Dgm: A deep learning algorithm for solving partial differential equations, Journal of computational physics 375 (2018) 1339--1364.

\bibitem{samaniego2020energy}
E.~Samaniego, C.~Anitescu, S.~Goswami, V.~M. Nguyen-Thanh, H.~Guo, K.~Hamdia, X.~Zhuang, T.~Rabczuk, An energy approach to the solution of partial differential equations in computational mechanics via machine learning: Concepts, implementation and applications, Computer Methods in Applied Mechanics and Engineering 362 (2020) 112790.

\bibitem{nguyen2021parametric}
V.~M. Nguyen-Thanh, C.~Anitescu, N.~Alajlan, T.~Rabczuk, X.~Zhuang, Parametric deep energy approach for elasticity accounting for strain gradient effects, Computer Methods in Applied Mechanics and Engineering 386 (2021) 114096.

\bibitem{chiu2022can}
P.-H. Chiu, J.~C. Wong, C.~Ooi, M.~H. Dao, Y.-S. Ong, Can-pinn: A fast physics-informed neural network based on coupled-automatic--numerical differentiation method, Computer Methods in Applied Mechanics and Engineering 395 (2022) 114909.

\bibitem{kharazmi2021hp}
E.~Kharazmi, Z.~Zhang, G.~E. Karniadakis, hp-vpinns: Variational physics-informed neural networks with domain decomposition, Computer Methods in Applied Mechanics and Engineering 374 (2021) 113547.

\bibitem{gao2021phygeonet}
H.~Gao, L.~Sun, J.-X. Wang, Phygeonet: physics-informed geometry-adaptive convolutional neural networks for solving parameterized steady-state pdes on irregular domain, Journal of Computational Physics 428 (2021) 110079.

\bibitem{sharma2022accelerated}
R.~Sharma, V.~Shankar, Accelerated training of physics informed neural networks (pinns) using meshless discretizations, arXiv preprint arXiv:2205.09332 (2022).

\bibitem{chakraborty2020simulation}
S.~Chakraborty, Simulation free reliability analysis: A physics-informed deep learning based approach, arXiv preprint arXiv:2005.01302 (2020).

\bibitem{zhang2022simulation}
C.~Zhang, A.~Shafieezadeh, Simulation-free reliability analysis with active learning and physics-informed neural network, Reliability Engineering \& System Safety 226 (2022) 108716.

\bibitem{jiang2024fourier}
Z.~Jiang, M.~Zhu, L.~Lu, Fourier-mionet: Fourier-enhanced multiple-input neural operators for multiphase modeling of geological carbon sequestration, Reliability Engineering \& System Safety (2024) 110392.

\bibitem{lu2019deeponet}
L.~Lu, P.~Jin, G.~E. Karniadakis, Deeponet: Learning nonlinear operators for identifying differential equations based on the universal approximation theorem of operators, arXiv preprint arXiv:1910.03193 (2019).

\bibitem{lu2021learning}
L.~Lu, P.~Jin, G.~Pang, Z.~Zhang, G.~E. Karniadakis, Learning nonlinear operators via deeponet based on the universal approximation theorem of operators, Nature Machine Intelligence 3~(3) (2021) 218--229.

\bibitem{li2020multipole}
Z.~Li, N.~Kovachki, K.~Azizzadenesheli, B.~Liu, A.~Stuart, K.~Bhattacharya, A.~Anandkumar, Multipole graph neural operator for parametric partial differential equations, Advances in Neural Information Processing Systems 33 (2020) 6755--6766.

\bibitem{xia2023maintenance}
L.~Xia, Y.~Liang, J.~Leng, P.~Zheng, Maintenance planning recommendation of complex industrial equipment based on knowledge graph and graph neural network, Reliability Engineering \& System Safety 232 (2023) 109068.

\bibitem{kovachki2021universal}
N.~Kovachki, S.~Lanthaler, S.~Mishra, On universal approximation and error bounds for fourier neural operators, Journal of Machine Learning Research 22~(290) (2021) 1--76.

\bibitem{tripura2023wavelet1}
T.~Tripura, A.~Awasthi, S.~Roy, S.~Chakraborty, A wavelet neural operator based elastography for localization and quantification of tumors, Computer Methods and Programs in Biomedicine 232 (2023) 107436.

\bibitem{tripura2023elastography}
T.~Tripura, A.~Awasthi, S.~Roy, S.~Chakraborty, A wavelet neural operator based elastography for localization and quantification of tumors, Computer Methods and Programs in Biomedicine 232 (2023) 107436.

\bibitem{rani2024generative}
J.~Rani, T.~Tripura, H.~Kodamana, S.~Chakraborty, Generative adversarial wavelet neural operator: Application to fault detection and isolation of multivariate time series data, arXiv preprint arXiv:2401.04004 (2024).

\bibitem{cao2023lno}
Q.~Cao, S.~Goswami, G.~E. Karniadakis, Lno: Laplace neural operator for solving differential equations, arXiv preprint arXiv:2303.10528 (2023).

\bibitem{tripura2023foundational}
T.~Tripura, S.~Chakraborty, A foundational neural operator that continuously learns without forgetting, arXiv preprint arXiv:2310.18885 (2023).

\bibitem{navaneeth2024physics}
N.~Navaneeth, T.~Tripura, S.~Chakraborty, Physics informed wno, Computer Methods in Applied Mechanics and Engineering 418 (2024) 116546.

\bibitem{daubechies1992ten}
I.~Daubechies, Ten lectures on wavelets, SIAM, 1992.

\bibitem{cotter2020uses}
F.~Cotter, Uses of complex wavelets in deep convolutional neural networks, Ph.D. thesis, University of Cambridge (2020).

\bibitem{selesnick2005dual}
I.~W. Selesnick, R.~G. Baraniuk, N.~C. Kingsbury, The dual-tree complex wavelet transform, IEEE signal processing magazine 22~(6) (2005) 123--151.

\bibitem{soin2024generative}
H.~Soin, T.~Tripura, S.~Chakraborty, Generative flow induced neural architecture search: Towards discovering optimal architecture in wavelet neural operator, arXiv preprint arXiv:2405.06910 (2024).

\bibitem{navaneeth2023stochastic}
N.~Navaneeth, S.~Chakraborty, Stochastic projection based approach for gradient free physics informed learning, Computer Methods in Applied Mechanics and Engineering 406 (2023) 115842.

\bibitem{mckean1970nagumo}
H.~McKean~Jr, Nagumo's equation, Advances in mathematics 4~(3) (1970) 209--223.

\bibitem{laing2009stochastic}
C.~Laing, G.~J. Lord, Stochastic methods in neuroscience, OUP Oxford, 2009.

\bibitem{lord2014introduction}
G.~J. Lord, C.~E. Powell, T.~Shardlow, An introduction to computational stochastic PDEs, Vol.~50, Cambridge University Press, 2014.

\bibitem{dejam2017pre}
M.~Dejam, H.~Hassanzadeh, Z.~Chen, Pre-darcy flow in porous media, Water Resources Research 53~(10) (2017) 8187--8210.

\bibitem{hussain2019approximate}
S.~Hussain, A.~Shah, S.~Ayub, A.~Ullah, An approximate analytical solution of the allen-cahn equation using homotopy perturbation method and homotopy analysis method, Heliyon 5~(12) (2019).

\bibitem{ma2017numerical}
L.~Ma, R.~Chen, X.~Yang, H.~Zhang, Numerical approximations for allen-cahn type phase field model of two-phase incompressible fluids with moving contact lines, Communications in Computational Physics 21~(3) (2017) 867--889.
\end{thebibliography}

\end{document}